\documentclass[10pt]{amsart}
\usepackage[dvipdfmx]{graphicx,color}
\usepackage{graphicx}

\newtheorem{theorem}{Theorem}[section]
\newtheorem{lemma}[theorem]{Lemma}

\theoremstyle{definition}
\newtheorem{definition}[theorem]{Definition}

\theoremstyle{remark}

\numberwithin{equation}{section}

\subjclass[2000]{Primary~47A75, Secondary~47A40}
\keywords{Quantum walk, Eigenvalue, Generalized eigenfunction, S-matrix}
\thanks{This work is supported by 	
Grant-in-Aid for Young Scientists (B) 16K17630, Japan Society for the Promotion of Science. }

\title[Generalized eigenfunctions and S-matrices for QW]
{Generalized eigenfunctions and scattering matrices for position-dependent quantum walks}

\author[H. Morioka]{Hisashi MORIOKA}
\address[H. Morioka]{Faculty of Science and Engineering,
Doshisha University, Tataramiyakodani 1-3, Kyotanabe, Kyoto, 610-0394, Japan}
\email{hisashimorioka@gmail.com}

\date{\today}

\begin{document}
\baselineskip 14pt
\maketitle

\begin{abstract}
We study the spectral analysis and the scattering theory for time evolution operators of position-dependent quantum walks.
Our main purpose of this paper is construction of generalized eigenfunctions of the time evolution operator.
Roughly speaking, the generalized eigenfunctions are not square summable but belong to $\ell^{\infty}$-space on ${\bf Z}$.
Moreover, we derive a characterization of the set of generalized eigenfunctions in view of the time-harmonic scattering theory.
Thus we show that the S-matrix associated with the quantum walk appears in the singularity expansion of generalized eigenfunctions.
\end{abstract}


\section{Introduction}
\subsection{Background}
Discrete time quantum walks (DTQWs or QWs for short) have been studied in various contexts of researches (see reviews of  DTQWs \cite{Me}, \cite{Am}, \cite{Kemp}, \cite{Kend}, \cite{Ko}, \cite{Ven}, \cite{Po}, \cite{Man} and references therein). 
Gudder \cite{Gu}, Meyer \cite{Me}, and Ambainis et al. \cite{ABNVW} considered one-dimensional DTQWs which are quantum versions of the random walk.
Nowadays, studies of DTQWs are flourishing in the context of quantum search algorithms or quantum computing  (see \cite{ABNVW} and Aharonov et al. \cite{AAKV}).

On the other hand, there is an abundance of recent works of position-dependent QWs (\cite{KLS}, \cite{EK1}, \cite{EK2}) in view of the spectral analysis and the quantum scattering theory.
In this point of view, we mention Asch et al. \cite{ABJ} and Suzuki \cite{Su}.
Asch et al. \cite{ABJ} studied Mourre's commutator method for unitary operators and proved the absence of the singular continuous spectrum.
Suzuki \cite{Su} proved that the Heisenberg operator of the position operator converges to the asymptotic velocity operator.
Richard et al. \cite{RST}, \cite{RST2} developed this direction of studies to anisotropic quantum walks, by using a two-Hilbert-space setting in order to obtain a weak limit theorem for quantum walks.
For commutator methods, see also Mourre \cite{Mor1} or Amrein et al. \cite{ABG}.
As has been shown by Cantero et al. \cite{CGM}, Segawa-Suzuki \cite{SS} and \cite{Su}, if the initial state has an overlap with an eigenspace of the time evolution operator, the associated QW has a localization. 
Some examples of localization with one-defect models are in \cite{CGM}, Konno et al. \cite{KLS} and Fuda et al. \cite{FuFuSu}.
Recently, Matsue et al. \cite{MMOS} constructed an example of counterparts of the resonant-tunneling effect for the quantum scattering with the double barrier potential.  
Morioka-Segawa \cite{MoSe} suggested a detection method of a kind of defects using embedded eigenvalues of time evolution operators.

In the time-harmonic scattering theory of self-adjoint operators like Schr\"{o}dinger operators, it is well-known that the wave operators are represented by distorted Fourier transformations which are defined by the spectral decompositions of Hamiltonians (see e.g. Yafaev \cite{Ya}).
Moreover, the scattering operator is defined through the wave operators and its Fourier transform is a direct integral of \textit{scattering matrices (S-matrices)}.
If we consider the Schr\"{o}dinger equation 
$$
(-\Delta +V -\lambda ) u=0 \quad \text{on} \quad {\bf R} , \quad \lambda >0,
$$
with $\mathrm{supp} \, V \subset [-R,R]$ for $R>0$, it is well-known that the wave function $ u$ is given by the generalized eigenfunction (which is not in $L^2$) of the form 
\begin{equation*}
u_{\pm} (x)=  \left\{
\begin{split}
T_{\pm} (\lambda ) e^{\pm i\sqrt{\lambda}x} &, \quad \pm x\geq  R, \\
e^{\pm i\sqrt{\lambda}x} + R_{\pm} (\lambda ) e^{\mp i \sqrt{\lambda}x} &, \quad \pm x\leq -  R.
\end{split}
\right.
\end{equation*}
The above constants $T_{\pm} (\lambda )$ and $R_{\pm} (\lambda )$ are called transmission coefficient and reflection coefficient. 
These values correspond the probability of transmission and reflection, respectively.
For one-dimensional Schr\"{o}dinger equation, the S-matrix is defined by $T_{\pm} (\lambda )$ and $R_{\pm} (\lambda )$.
More generally, if we consider $H_0 = -\Delta $, $H=H_0 +V$, $\lim _{|x| \to \infty} \langle x \rangle^{\rho} V(x)=0$ for some constants $\rho >1$ on ${\bf R}^d$, the wave operator is defined by
$$
W_{\pm} = {\mathop{{\rm s-lim}}_{t\to\pm \infty}} \, e^{itH} e^{-itH_0} \quad \text{in} \quad L^2 ({\bf R}^d ) ,
$$
and then the scattering operator is given by
$$
{\bf S} = (W_+ )^* W_- .
$$
Physically, the operator ${\bf S}$ relates the behavior of the quantum particle at $t\to -\infty $ and at $t\to \infty $.
The S-matrix is given by the Fourier transform $\widehat{{\bf S}}$ of ${\bf S}$ of the form
$$
\widehat{{\bf S}} = \int _0^{\infty} \oplus \widehat{{\bf S}} (\lambda ) d\lambda ,
$$
where $\widehat{{\bf S}} (\lambda )$ is a unitary operator on $L^2 (S^{d-1} )$.
Moreover, it is well-known that the S-matrix is represented explicitly by the distorted Fourier transformation associated with the Hamiltonian $H$.

\subsection{Model of DTQW}
In this paper, we consider position-dependent quantum walks on ${\bf Z} $.
Let $ \mathcal{H} = \ell^2 ({\bf Z} ; {\bf C}^2 )$ be the space of states.
We define the unitary operator $U$ by 
\begin{equation}
(U\psi )(x)= C_1 (x+1) \psi (x+1) + C_2 (x-1) \psi (x-1) , \quad x\in {\bf Z} ,
\label{S1_def_QW}
\end{equation}
where $\psi \in \mathcal{H} $ and 
$$
C_1 (x)= \left[ \begin{array}{cc}
a(x) & b(x) \\ 0 & 0 \end{array} \right] , \quad C_2 (x)= \left[ \begin{array}{cc}
 0 & 0 \\ c(x) & d(x) \end{array} \right] .
$$
We assume $ C(x) := C_1 (x) + C_2 (x) \in U(2)$ for every $x\in {\bf Z} $.
We denote by $C$ the operator of multiplication by $C(x)$ for each $x\in {\bf Z}$ i.e. $(C\psi ) (x)= C(x) \psi (x) $ for $\psi \in \mathcal{H}$. 
We call $C$ the \textit{coin operator}.  
The operator $U$ is written by $ U=SC$ where $S$ is the shift operator defined by 
$$
(S\psi )(x)= \left[ \begin{array}{c}
\psi ^{(0)} (x+1) \\ \psi ^{(1)} (x-1) \end{array} \right] .
$$
The operator $U$ is a time evolution operator for quantum walks.
In fact, taking an initial state $ \psi_0 \in \mathcal{H} $, the state at time $t\in {\bf Z}$ is given by $\psi (t,\cdot )= U^t \psi _0 $.
We call this time evolution \textit{position-dependent quantum walk}.

We consider the operator $U$ in view of the scattering theory.
Thus we introduce the corresponding position-independent quantum walk.
Let $ U_0 = S C_0 $ be a unitary operator on $\mathcal{H} $ where 
$$
C_0 = \left[ \begin{array}{cc}
a_0 & b_0 \\ c_0 & d_0 \end{array} \right] \in U(2) .
$$
Moreover, we use the following representation of $C_0$ introduced in \cite{RST} : 
\begin{equation}
C_0 = e^{i\gamma /2} \left[ \begin{array}{cc}
p e^{i ( \alpha - \gamma /2 )} & q e^{i(\beta - \gamma /2)} \\ -q e^{-i (\beta - \gamma /2 )} & p e^{-i(\alpha -\gamma /2 )} \end{array} \right] ,
\label{S1_eq_C0}
\end{equation}
where $ p \in (0,1]$, $q\in [0,1)$ with $p^2 + q^2 =1$ and $ \alpha , \beta , \gamma \in {\bf R} /(2\pi {\bf Z} )$.
Throughout of the paper, we adopt the following assumptions : 

\medskip 

{\bf (A-1)} There exist constants $ \epsilon_0 , M >0$ such that 
$$
\| C(x) -C_0 \| _{\infty} \leq M e^{-\epsilon_0 \langle x \rangle} , \quad x\in {\bf Z} , \quad \langle x \rangle = \sqrt{1+x^2 } .
$$

{\bf (A-2)} There exists a constant $\delta >0$ such that $ |a(x)| \geq \delta $ for all $x\in {\bf Z} $.

\medskip

Here $\| \cdot \|_{\infty}$ is the norm of $2\times 2$ matrices defined by 
$$
\| A\| _{\infty} = \max _{1\leq j,k \leq 2} | a_{jk} | , \quad A=[ a_{jk} ] _{1\leq j,k \leq 2} .
$$

\subsection{Purpose, summary of result and plan of the paper} 
In this paper, we study the spectral analysis and the scattering theory for $U$ and $U_0$.
In particular, we construct the S-matrix associated with the wave operator 
$$
W_{\pm} = {\mathop{{\rm s-lim}}_{t\to\pm \infty}} \, U^{-t} U_0^t  \quad \text{in} \quad \mathcal{H} .
$$
In order to do this, we adopt the time-harmonic approach; that is the limiting absorption principle and the construction of the distorted Fourier transformation in view of the spectral decomposition of $U$.
As is well-known, the spectral analysis which will be done in this paper usually acts for self-adjoint Hamiltonians.
Since $U$ and $U_0$ are unitary operators on $\mathcal{H}$, there exist self-adjoint operators $H$ and $H_0$ on $\mathcal{H}$ such that $U = e^{iH} $ and $U_0 =e^{iH_0}$. 
However, $H$ and $H_0$ are non-local operators in general and they  are not easy to analyze in view of the spectral theory, although its concrete formulas have been given.
For this topic for position-independent cases, see Tate \cite{Ta}.
We also mention Segawa-Suzuki \cite{SS}, Fuda et al. \cite{FuFuSu} and \cite{FuFuSu2}.
The authors used the discriminant operator which is a kind of Hamiltonians generating the time evolution operators of some DTQWs.
Fortunately, the determinant of the Fourier multiplier of $U_0$ is similar to that of discrete Laplacian.
Then we shall study $U$ and $U_0$ directly, without passing to their generating Hamiltonians $H$ and $H_0$. 

The limiting absorption principle and the distorted Fourier transformation enable us to derive the generalized eigenfunction of $U$ with the continuous spectrum.
Roughly speaking, the generalized eigenfunction is not in $\mathcal{H}$ but in $ \ell^{\infty} ({\bf Z} ; {\bf C}^2 )$.
This is a generalization of tunneling solutions of DTQWs given in \cite{MMOS}.
The S-matrix will be represented by the distorted Fourier transformation and the generalized eigenfunction.

The plan of this paper is as follows.
In \S 2, we introduce functional spaces and recall some fundamental materials of spectral properties.
In \S 3, we prove the limiting absorption principle for $U_0$ passing to the Fourier series.
The main analytical tool of the proof is a division theorem for distribution on the torus.
In order to discuss it, we use complex contour integrations. 
Some basic formulas are gathered in the appendix.
Moreover, the generalized eigenfunction for $U_0$ is defined here. 
We will prove the uniqueness of a kind of solutions to the equation $(U_0 -e^{i\theta} )u=f$ and the characterization of the set of solutions to the equation $( U_0 - e^{i\theta} )u=0$ in the sense of Agmon-H\"{o}rmander's $ \mathcal{B} $-$\mathcal{B}^* $ space (\cite{AgHo}).
In \S 4, we prove the limiting absorption principle and construct the generalized eigenfunction for $U$ based on the argument in \S 3. 
The generalized eigenfunction is defined as the adjoint operator of the distorted Fourier transformation associated with $U $.
The main result of the paper is derived in \S 5.
Based on the definition of the wave operators by Suzuki \cite{Su} and Richard et al. \cite{RST}, \cite{RST2}, we define the scattering operator and associated S-matrix.
We also derive a concrete formula of the S-matrix by using the distorted Fourier transformation.
This formula relates the singularity expansion of the generalized eigenfunction and the S-matrix.
Our main theorems are Theorems \ref{S5_thm_Smatrix} and \ref{S5_thm_rangeF+}.

\subsection{Remark}
Except for Theorems \ref{S2_thm_embev} and \ref{S4_thm_absenceofef}, our arguments work under an assumption weaker than (A-1) and (A-2).
In fact, as has been mentioned above, Suzuki \cite{Su} and Richard et al. \cite{RST}, \cite{RST2} assumed the short-range condition i.e. the right-hand side of the condition (A-1) can be replaced by $M\langle x \rangle ^{-1-\epsilon_0 } $ for some constants $\epsilon_0 , M >0$. 
We should mention that Maeda et al. \cite{MSSS} also proved some parts of this paper under the short-range condition in order to prove a dispersive estimate for DTQWs. 
However, our aim which has been mentioned above and the method which will be used in this paper are different from \cite{MSSS}.
The pair of Theorems \ref{S2_thm_embev} and \ref{S4_thm_absenceofef} is an analogue of Rellich's uniqueness theorem for the Helmholtz type equations on ${\bf R}^d $ (\cite{Re}, \cite{Vek}, \cite{Tr}, \cite{Li1}, \cite{Li2}, \cite{Ho}, \cite{Mu}, \cite{RaTa}).
Note that the Rellich's uniqueness theorem for the discrete Schr\"{o}dinger operator is still open problem if we assume that the perturbation is polynomially decaying.
If the perturbation is finite rank or exponentially decaying, Rellich's theorem holds for discrete Schr\"{o}dinger operators (\cite{IsMo}, \cite{Ves}, \cite{AIM}).

\subsection{Notation}
The notation in this paper is as follows.
${\bf T} := {\bf R} / (2\pi {\bf Z} )$ denotes the torus.
For $f\in \mathcal{S}' ({\bf R})$, $\widetilde{f} (\xi )$ denotes its Fourier transform
$$
\widetilde{f} (s )= (2\pi)^{-1/2} \int _{{\bf R}} e^{-is \xi} f(\xi ) d\xi , \quad s\in {\bf R} .
$$
For $f =\{ f(x)\} _{x\in {\bf Z}}$, we define the mapping $\mathcal{U} : \ell^2 ({\bf Z} ) \to L^2 ({\bf T} )$ by
$$
\widehat{f} (\xi) :=(\mathcal{U} f) (\xi ) = (2\pi )^{-1/2} \sum _{x\in {\bf Z}} e^{-ix\xi } f (x) ,\quad \xi \in {\bf T} .
$$
However, we also use the notation $ \mathcal{U} $ as the Fourier series such that $f\not\in \mathcal{H}$.
For $\widehat{g} \in \mathcal{S}' ({\bf T} )$, $g (x)= (\mathcal{U}^* \widehat{g})(x)$ for $x\in {\bf Z}$ is its Fourier coefficient
$$
g(x)= (2\pi)^{-1/2} \int _{{\bf T}} e^{ix\xi} \widehat{g} (\xi ) d\xi .
$$
For a vector $ {\bf v} = {}^t (v^{(0)} , v^{(1)} ) \in {\bf C}^2$, we define the norm $| {\bf v} | _{{\bf C}^2} = \sqrt{ |v^{(0)} |^2 +|v^{(1) }|^2 }$.
For ${\bf u} , {\bf v} \in {\bf C}^2$, we denote by $({\bf u}, {\bf v} )_{{\bf C}^2} $ the standard inner product of ${\bf C}^2 $.
In this paper, we often consider ${\bf C}^2$-valued functions $f= {}^t (f^{(0)} , f^{(1)} )$ on ${\bf Z}$, ${\bf T}$ or ${\bf R}$.
We also use the notations $\widehat{f} = {}^t (\widehat{f}^{(0)} , \widehat{f} ^{(1)} ) $ and $
\widetilde{f} = {}^t (\widetilde{f}^{(0)} , \widetilde{f} ^{(1)} )$. 


\section{Preliminaries}

\subsection{Sobolev and Besov spaces}

We introduce functional spaces for ${\bf C}^2$-valued functions.
For ${\bf C}$-valued functions, these functional spaces can be defined by the similar way, replacing $|\cdot | _{{\bf C}^2} $ by $|\cdot |$.

Let $ r_{-1} =0$, $r_j = 2^j $ for $j\geq 0$.
We define the space $ \mathcal{B} ({\bf R})$ to be the set of functions $f= {}^t (f^{(0)} , f^{(1)}) $ on ${\bf R}$ having the norm
$$
\| f\| _{\mathcal{B}({\bf R})} = \sum _{j=0}^{\infty} r_j^{1/2} \left( \int _{\Omega_j } | \widetilde{f} (s)|^2_{{\bf C}^2} ds \right)^{1/2} ,
$$
where $\Omega_j = \{ s\in {\bf R} \ ; \ r_{j-1} \leq |s | < r_j \} $.
The (equivalent) norm of its dual space $ \mathcal{B}^* ({\bf R} )$ is 
$$
\| u\| _{\mathcal{B}^* ({\bf R} )} ^2 = \sup _{R>1} \frac{1}{R} \int_{|s|<R} |\widetilde{u} (s)|^2 _{{\bf C}^2} ds .
$$
The space $ \mathcal{B} _0^* ({\bf R} )$ is defined by
$$
\mathcal{B} _0^* ({\bf R} ) = \left\{ u\in \mathcal{B}^* ({\bf R} ) \ ; \ \lim_{R\to \infty } \frac{1}{R} \int _{|s|<R} |\widetilde{u} (s) |^2_{{\bf C}^2} ds =0 \right\} .
$$
The Sobolev space $H^{\sigma} ({\bf R} )$ is defined by
$$
H^{\sigma} ({\bf R} )=\left\{ u\in \mathcal{S}' ({\bf R} ) \ ; \ \| (1+s^2 )^{\sigma /2} \widetilde{u} (s) \|_{L^2 ({\bf R} ; {\bf C}^2 )} < \infty \right\} , \quad \sigma \in {\bf R} .
$$

We can define $ \mathcal{B} ({\bf T})$, $\mathcal{B}^* ({\bf T})$ and $H^{\sigma} ({\bf T} )$, taking a partition of unity.
Let $\{ \chi_j \} $ be a partition of unity on ${\bf T} $ where the support of $\chi_j $ is sufficiently small.
Then we define
$$
\| \widehat{f} \| _{\mathcal{B} ({\bf T} )} = \sum _{j} \| \chi_j \widehat{f} \| _{\mathcal{B} ({\bf R} )} , \quad \| \widehat{u} \| _{\mathcal{B}^* ({\bf T} )} = \sum _{j} \| \chi_j \widehat{u} \| _{\mathcal{B}^* ({\bf R} )} .
$$
We define $\widehat{u} \in \mathcal{B}_0^* ({\bf T})$ by $\chi_j \widehat{u} \in \mathcal{B}_0^* ({\bf R} )$ for all $j$.
The Sobolev space $H^{\sigma} ({\bf T} )$ is defined by the similar way.

The functional spaces $ \mathcal{B} ({\bf Z} )$, $ \mathcal{B}^* ({\bf Z} )$ are defined by the norms
\begin{gather*}
\| f \| _{\mathcal{B} ({\bf Z})} = \sum _{j=0}^{\infty} r_j^{1/2} \left( \sum _{r_{j-1} \leq |x| < r_j} | f(x) |_{{\bf C}^2}^2 \right)^{1/2} , \\
\| u\| _{\mathcal{B}^* ({\bf Z} )}^2 = \sup _{R>1} \frac{1}{R} \sum _{|x| <R} | u(x)|^2_{{\bf C}^2} .
\end{gather*}
Thus $\mathcal{B}_0^* ({\bf Z} )$ is defined by 
$$
\mathcal{B}_0^* ({\bf Z} )= \left\{ u\in \mathcal{B}^* ({\bf Z}) \ ; \ \lim _{R\to \infty} \frac{1}{R} \sum _{|x|<R} | u(x)|^2_{{\bf C}^2} =0 \right\} .
$$
The weighted $\ell^2$-spaces $ \ell^{2,\sigma} ({\bf Z} )$ are given by the set of $u=\{ u(x) \}_{x\in {\bf Z}}$ satisfying
$$
\| u\| ^2_{\ell^{2,\sigma} ({\bf Z} )} = \sum _{x\in {\bf Z}} (1+|x|^2 )^{\sigma} |u(x)|^2 _{{\bf C}^2} <\infty , \quad \sigma \in {\bf R}  .
$$

Using the Fourier series, we can relate these functional spaces on ${\bf Z}$ and on ${\bf T}$.
We define the position operator $X$ by
$$
(Xf)(x)= x f(x) , \quad x\in {\bf Z} .
$$
Then we have 
$$
\widehat{X} = \mathcal{U} X \mathcal{U}^* = i \frac{d}{d\xi} , \quad \widehat{X}^2 = - \frac{d^2}{d\xi^2} ,
$$
with the periodic boundary condition on ${\bf T}$.
For a self-adjoint operator $A$, $ \chi ( a\leq A<b)$ denotes the operator $\chi_I ( A )$ where $ \chi _I (\lambda )$ is the characteristic function of the interval $I=[a,b)$.
The operators $\chi ( A<a)$ and $ \chi (A \geq b)$ are defined by the similar way.
Let us introduce (equivalent) norms :
\begin{gather*}
\| \widehat{f} \| _{\mathcal{B} ({\bf T} )} = \sum _{j=0}^{\infty} r_j^{1/2} \| \chi ( r_{j-1} \leq |\widehat{X} | < r_j ) \widehat{f} \| _{L^2 ({\bf T} ; {\bf C}^2 )} , \\
\| \widehat{u} \| ^2 _{\mathcal{B}^* ({\bf T})} = \sup _{R>1} \frac{1}{R} \| \chi (|\widehat{X} |<R) \widehat{u} \| ^2 _{L^2 ({\bf T} ;{\bf C}^2 )} ,\\
\| \widehat{u} \| _{H^{\sigma} ({\bf T} )} = \| (1+\widehat{
X}^2  )^{\sigma /2} \widehat{u} \| _{L^2 ({\bf T} ;{\bf C}^2 )} , \quad \sigma \in {\bf R}  .
\end{gather*}
The space $ \mathcal{B}_0^* ({\bf T})$ is rewritten as 
$$
\mathcal{B}_0^* ({\bf T}) = \left\{ u\in \mathcal{B}^* ({\bf T}) \ ; \ \lim_{R\to \infty} \frac{1}{R} \| \chi (|\widehat{X} |<R) \widehat{u} \| ^2 _{L^2 ({\bf T} ;{\bf C}^2 )} =0 \right\} .
$$
Then we have 
\begin{gather*}
f\in \mathcal{B} ({\bf Z}) \Longleftrightarrow \widehat{f} \in \mathcal{B} ({\bf T}) ,\\
u\in \ell^{2,\sigma} ({\bf Z}) \Longleftrightarrow \widehat{u} \in H^{\sigma} ({\bf T} ) ,\\
u\in \mathcal{B}^* ({\bf Z}) \Longleftrightarrow \widehat{u} \in \mathcal{B}^* ({\bf T}) ,\\
u\in \mathcal{B}_0^* ({\bf Z}) \Longleftrightarrow \widehat{u} \in \mathcal{B}^*_0 ({\bf T}) .
\end{gather*}
In particular, we have the following inclusion relation.
For details, see \cite{AgHo}, Chapter 14 in \cite{Zw} and \S 3 in \cite{IsMo}.

\begin{lemma}
For $\sigma >1/2 $, we have
\begin{gather*}
\ell^{2,\sigma} ({\bf Z}) \subset \mathcal{B} ({\bf Z}) \subset \ell^{2,1/2} ({\bf Z}) \subset \ell^2 ({\bf Z} )\subset \ell^{2,-1/2} ({\bf Z} ) \subset \mathcal{B}^* ({\bf Z}) \subset \ell^{2,-\sigma} ({\bf Z}) , \\
H^{\sigma} ({\bf T} ) \subset \mathcal{B} ({\bf T}) \subset H^{1/2} ({\bf T}) \subset L^2 ({\bf T}) \subset H^{-1/2} ({\bf T}) \subset \mathcal{B}^* ({\bf T}) \subset H^{-\sigma} ({\bf T}).
\end{gather*}

\label{S2_lem_inclusion}
\end{lemma}

For $u\in \mathcal{B}^* ({\bf Z})$ and $f\in \mathcal{B} ({\bf Z})$, $(u,f)$ denotes the coupling between $u$ and $f$ :
$$
(u,f)= \sum _{x\in {\bf Z}} (u(x) , f(x)) _{{\bf C}^2} .
$$
For $ \mathcal{B} ({\bf T})$ and $\mathcal{B} ^* ({\bf T})$, or $\mathcal{B} ({\bf R})$ and $\mathcal{B}^* ({\bf R})$, we use the same notation in order to represent the similar coupling.

Finally, we introduce basic inequalities for $\mathcal{B}$.

\begin{lemma}
(1) For $f\in \mathcal{B} ({\bf R})$, we have
$$
\int _{-\infty} ^{\infty} | \widetilde{f} (s) |_{{\bf C}^2} ds \leq \sqrt{2} \| f\| _{\mathcal{B} ({\bf R} )} .
$$
(2) Let $f\in \mathcal{B} ({\bf Z})$.
For any $\xi \in {\bf T}$, $ \widehat{f}$ satisfies
$$
|\widehat{f} (\xi )|\leq \frac{1}{\sqrt{\pi}} \| f\| _{\mathcal{B} ({\bf Z})} .
$$ 

\label{S2_lem_ell1}
\end{lemma}

Proof.
It follows from the Cauchy-Schwartz inequality that
\begin{gather*}
\begin{split}
\int _{-\infty}^{\infty} | \widetilde{f} (s) | _{{\bf C}^2} ds & = \sum _{j\geq 0} \int _{\Omega_j} |\widetilde{f} (s) |_{{\bf C}^2} ds \\
& \leq \sqrt{2} \sum _{j\geq 0} r_j^{1/2} \left( \int _{\Omega_j} | \widetilde{f} (s) |^2 _{{\bf C}^2} ds \right)^{1/2} .
\end{split}
\end{gather*}
The right-hand side is equal to $\sqrt{2} \| f\| _{\mathcal{B} ({\bf R})} $. 
This proves the assertion (1).

For $f\in \mathcal{B} ({\bf Z})$, we have 
$$
|\widehat{f} (\xi )| \leq \frac{1}{\sqrt{2\pi}} \sum_{x\in {\bf Z}} |f(x)| \leq \frac{1}{\sqrt{\pi}} \| f\| _{\mathcal{B} ({\bf Z})} , \quad \xi \in {\bf T} ,
$$
by the similar way.
Thus we also have the assertion (2).
\qed

\medskip

Throughput of this paper, we often discuss in $ \mathcal{B}$-$\mathcal{B}^*$ spaces i.e. outside the space of states $ \mathcal{H} $.
Then we naturally extend $U$ and $U_0$ to $ \mathcal{B}^*$ space.
The operators $U$ and $U_0$ are not unitary on $\mathcal{B}^*$, but the notation $U^*$ and $  U_0^*$ will be used when there is no afraid of confusion, in the sense of 
$$
(U^* u)(x)= C(x)^* \left[ \begin{array}{c} u^{(0)} (x-1) \\ u^{(1)} (x+1) \end{array} \right] , \quad x\in {\bf Z},
$$
for $u\in \mathcal{B}^* ({\bf Z})$.
The operator $U_0^*$ on $\mathcal{B}^* ({\bf Z})$ is given by the same manner.


\subsection{Essential spectrum}

Spectral decompositions of unitary operators are discussed in \cite{RST} and \cite{MoSe}.
Here we recall some basic properties of spectra of unitary operators.

Let $ U $ be a unitary operator on a Hilbert space $ \mathcal{H} $.
We denote by $\sigma (U)$ the spectrum of $U$.
There exists a spectral decomposition $E_U ( \theta )$ for $ \theta \in {\bf R} $ such that 
$$
U= \int_0^{2\pi} e^{i\theta} dE_U (\theta ) ,
$$
where $E_U (\theta )$ is extended to be zero for $ \theta \in ( - \infty , 0 )$ and to be $1$ for $ \theta \in [2\pi , \infty )$.
Thus the orthogonal decomposition of $ \mathcal{H} $ associated with $U$ is given by 
$$
\mathcal{H} = \mathcal{H}_p (U) \oplus \mathcal{H}_{sc} (U ) \oplus \mathcal{H}_{ac} (U) ,
$$
where $ \mathcal{H}_p (U) $, $ \mathcal{H}_{sc} (U) $ and $ \mathcal{H}_{ac} (U)$ are orthogonal projections on pure point, singular continuous and absolutely continuous subspaces of $\mathcal{H} $, respectively.
Note that $ \mathcal{H}_p (U)$ is the closure of all eigenspaces.
Then we define $\sigma_p (U)$ as the set of eigenvalues of $U$ and
$$
 \sigma_{sc} (U)= \sigma (U| _{\mathcal{H}_{sc} (U)} ) ,\quad \sigma_{ac} (U)= \sigma (U| _{\mathcal{H}_{ac} (U)} ),
$$
and we call them point spectrum, singular continuous spectrum and absolutely continuous spectrum, respectively.

We also define discrete spectrum and essential spectrum of $U$.
The discrete spectrum $\sigma_d (U)$ is the set of isolated eigenvalues of $U$ with finite multiplicities.
The essential spectrum $\sigma_{ess} (U)$ is defined by $ \sigma_{ess} (U)= \sigma (U) \setminus \sigma_d (U)$.
If $ \lambda \in \sigma_{ess} (U)$, $\lambda$ is either an eigenvalue with infinite multiplicity or an accumulation point of $ \sigma (U)$.
We state the following property.
For its proof, see Lemmas 2.1 and 2.2 in \cite{MoSe}.

\begin{lemma}
Let $U'$ and $U$ be unitary operators on $\mathcal{H}$.
If $U'-U$ is a compact operator on $\mathcal{H}$, we have $ \sigma_{ess} (U') = \sigma_{ess} (U)$. 
\label{S2_lem_singularseq}
\end{lemma}

Now we turn to the quantum walk.
Let $U=SC$ be define by (\ref{S1_def_QW}) and $U_0 =SC_0$ be given by $C_0$ in (\ref{S1_eq_C0}).
Letting 
$$
\widehat{U}_0 = \mathcal{U} U_0 \mathcal{U}^* ,
$$
it follows that $\widehat{U}_0$ is the operator of multiplication on ${\bf T}$ by the unitary matrix
$$
\widehat{U}_0 (\xi )= e^{i\gamma /2} \left[ \begin{array}{cc}
p e^{i(\alpha -\gamma /2)} e^{i\xi} & q e ^{i(\beta -\gamma /2)} e^{i\xi} \\ -q e^{-i (\beta - \gamma /2)} e^{-i\xi} & p e ^{-i (\alpha - \gamma /2)} e^{-i\xi} \end{array} \right] , \quad \xi \in {\bf T} .
$$
Then we have for any $ z \in {\bf C} $
\begin{gather}
\begin{split}
p(\xi , z ) := & \, \mathrm{det} (\widehat{U} _0 (\xi ) -e^{iz} ) \\
= & \, 2p e^{i(z+\gamma /2)} \left( -\cos \left( \xi + \alpha -\frac{\gamma}{2} \right) + \frac{1}{p} \cos \left( z-\frac{\gamma}{2}  \right)  \right) .
\end{split}
\label{S2_eq_determinant}
\end{gather}

\begin{lemma}
Let $ J_{\gamma} = J_{\gamma ,1} \cup J_{\gamma ,2}$ where 
\begin{gather*}
J_{\gamma ,1} = [ \arccos p + \gamma /2 , \pi - \arccos p + \gamma /2 ], \\
J_{\gamma ,2} = [ \pi + \arccos p + \gamma /2 , 2\pi - \arccos p + \gamma /2 ] .
\end{gather*}
Then we have $ \sigma (U_0 )= \sigma _{ac} (U_0) = \{ e^{i\theta} \ ; \ \theta \in J_{\gamma} \} $.
In particular, we have $\sigma (U_0)= S^1 $ if $p=1$.
\label{S2_lem_essU0}
\end{lemma}

Proof.
See Lemma 4.1 in \cite{RST}.
\qed

\medskip

In view of the assumptions (A-1) and (A-2), $U-U_0$ is compact on $\mathcal{H} $.
The following lemma follows from Lemma \ref{S2_lem_singularseq}.

\begin{lemma}
We have $\sigma_{ess} (U)= \sigma_{ess} (U_0)= \{ e^{i\theta} \ ; \ \theta \in J_{\gamma} \} $.
\label{S2_lem_essU}
\end{lemma}

Under the assumptions (A-1) and (A-2), we can also prove the absence of eigenvalues embedded in the continuous spectrum as follows.
 
\begin{theorem}
Let $ \mathcal{T} = \{ e^{i\theta} \ ; \ \theta \in J_{\gamma , \mathcal{T}} \} $ where 
$$
J_{\gamma , \mathcal{T}} = \left\{ \begin{split}
&\arccos p + \gamma /2 , \pi - \arccos p + \gamma /2 , \\
&\pi + \arccos p + \gamma /2 , 2\pi - \arccos p + \gamma /2 \end{split} \right\} .
$$
There is no eigenvalue of $U$ in $\sigma_{ess} (U) \setminus \mathcal{T} $.
\label{S2_thm_embev}
\end{theorem} 

Proof. See Theorem 1.2 in \cite{MoSe}. 
\qed

\medskip

The set $ \mathcal{T} $ consists of thresholds in $\sigma_{ess} (U)$.
Let 
\begin{gather}
M(\theta )=\{ \xi \in {\bf T} \ ; \ p(\xi ,\theta )=0 \} ,\label{S2_def_fermi} \\
M_{reg} (\theta )= \left\{ \xi \in {\bf T} \ ; \ p(\xi ,\theta )=0 , \ \frac{\partial p}{\partial \xi} (\xi ,\theta ) \not= 0 \right\} ,\label{S2_def_fermireg} \\
M_{sng} (\theta )= \left\{ \xi \in {\bf T} \ ; \ p(\xi ,\theta )=0 , \ \frac{\partial p}{\partial \xi} (\xi ,\theta ) = 0 \right\} .\label{S2_def_fermisng} 
\end{gather}

\begin{lemma}
If $ e^{i\theta} \in \sigma_{ess} (U) \setminus \mathcal{T}$, we have $ M(\theta )= M_{reg} (\theta )$.
On the other hand, if $ e^{i\theta} \in \mathcal{T}$, we have $ M(\theta )= M_{sng} (\theta )$.
\label{S2_lem_fermi}
\end{lemma}

Proof.
See Lemma 3.1 in \cite{MoSe}.
\qed



\section{Position-independent QW}

\subsection{Resolvent estimate}

Here we will derive a division theorem and its micro-local properties of the resolvent operator
\begin{equation}
R_0 (z)= ( U_0 -e^{iz} )^{-1} , \quad z \in {\bf C} \setminus {\bf R}  .
\label{S3_def_resolvent}
\end{equation}
Letting $ \widehat{R}_0 (z)= \mathcal{U} R_0 (z) \mathcal{U}^* $, $\widehat{R}_0 (z)$ is the operator of multiplication by the matrix 
\begin{equation}
\widehat{R} _0 (\xi, z) = \frac{1}{p(\xi ,z)} \left[ \begin{array}{cc} pe^{ -i ( \alpha -\gamma )} e^{-i\xi} -e^{iz }  & -q e^{i\beta} e^{i\xi } \\ q e^{-i (\beta -\gamma )} e^{-i\xi } & p e^{i\alpha} e^{i\xi} - e^{iz} \end{array} \right] ,
\label{S3_def_resolventtorus}
\end{equation}
for $\xi \in {\bf T} $.

We denote by $ \lambda_1 (\xi )$ and $\lambda_2 (\xi )$ the eigenvalues of the matrix $\widehat{U}_0 (\xi )$.
In view of (\ref{S2_eq_determinant}), we can see $ \lambda_1 (\xi )= e^{i\theta (\xi )} $ and $ \lambda_2 (\xi )= e^{i(\gamma - \theta (\xi ))} (= e^{ i(2\pi + \gamma -\theta (\xi ))} ) $ where 
\begin{equation}
\theta (\xi )= \frac{\gamma}{2} + \arccos \left( p \cos \left( \xi + \alpha -\frac{\gamma}{2} \right) \right) .
\label{S3_eq_thetaxi}
\end{equation}
Moreover, we have
\begin{equation}
J_{\gamma ,1} = \{ \theta (\xi ) \ ; \ \xi \in {\bf T} \}   , \quad  J_{\gamma ,2} = \{ 2\pi + \gamma -\theta (\xi ) \ ; \ \xi \in {\bf T} \}  ,
\label{S3_eq_range}
\end{equation}
and 
\begin{gather}
\frac{d \lambda_1}{d\xi } (\xi )= i\lambda_1 (\xi ) \frac{p\sin (\xi + \alpha -\gamma /2 ) }{ \sqrt{1-p^2 \cos^2 (\xi + \alpha -\gamma /2 ) }} , \label{S3_eq_dlambda1}  \\
\frac{d \lambda_2}{d\xi } (\xi )=- i\lambda_2 (\xi ) \frac{p\sin (\xi + \alpha -\gamma /2 ) }{ \sqrt{1-p^2 \cos^2 (\xi + \alpha -\gamma /2 ) }} .\label{S3_eq_dlambda2}
\end{gather}
Then if $ \xi $ varies in $ \cup _{\omega \in J} M(\omega )$ for a compact interval $J \subset J_{\gamma} \setminus J_{\gamma , \mathcal{T}} $, it follows that $\lambda_1 (\xi )$ and $\lambda_2 (\xi )$ are distinct, and $( d\lambda_1 / d\xi )(\xi )$ and $ (d\lambda_2 / d\xi )(\xi )$ do not vanish. 
For $ \theta \in J_{\gamma,j} \setminus J_{\gamma , \mathcal{T}} $, $M(\theta )$ consists of two distinct points on ${\bf T}$.

Let $ P_j (\xi )$ be the orthogonal projection on the eigenspace of $ \lambda_j (\xi )$. 
Thus $\widehat{U}_0 (\xi )$ can be represented as 
$$
\widehat{U}_0 (\xi )= \lambda_1 (\xi ) P_1 (\xi ) + \lambda_2 (\xi ) P_2 (\xi ) .
$$
$P_j (\xi )$ is smooth in a small neighborhood of each point in $ \cup _{\omega \in J} M(\omega )$.

Let $ \theta \in J_{\gamma} \setminus J_{\gamma , \mathcal{T} } $ and $\xi (\theta ) \in M(\theta )$. 
We take $\chi \in C^{\infty} ({\bf T} )$ such that $ \chi (\xi (\theta ))=1$ with small support.
For $ \widehat{u} = \widehat{R}_0 (z) \widehat{f} $ with $ \widehat{f} \in \mathcal{B} ({\bf T})$, we have 
\begin{equation}
( \lambda_j (\xi )-e^{iz} ) \chi (\xi ) P_j (\xi ) \widehat{u} (\xi )=\chi (\xi ) P_j (\xi ) \widehat{f} (\xi ) , \quad z\in {\bf C} .
\label{S3_eq_helmholtz}
\end{equation}

Obviously, we see the following lemma.

\begin{lemma}
If $ \theta \in J_{\gamma , 1} \setminus J_{\gamma , \mathcal{T}} $, $\lambda_2 (\xi )-e^{i\theta} $ is invertible. 
If $ \theta \in J_{\gamma , 2} \setminus J_{\gamma , \mathcal{T}} $, $\lambda_1 (\xi )-e^{i\theta} $ is invertible.
\label{S3_lem_distinct}
\end{lemma}

In view of (\ref{S3_eq_thetaxi}), we introduce a change of variable $\xi \mapsto \eta $ in a small neighborhood of a point $\xi (\theta )\in M(\theta )$ as follows.
For $ \theta \in J_{\gamma ,1} \setminus J_{\gamma , \mathcal{T} } $, let 
\begin{gather*}
\begin{split}
\eta &=  \arccos \left( p\cos \left( \xi + \alpha -\frac{\gamma}{2} \right)  \right) - \arccos \left( p\cos \left( \xi (\theta ) + \alpha -\frac{\gamma}{2} \right)  \right) \\
&= \arccos \left( p\cos \left( \xi + \alpha -\frac{\gamma}{2} \right)  \right) +\frac{\gamma}{2} - \theta .
\end{split}
\end{gather*}
For $ \theta \in J_{\gamma ,2} \setminus J_{\gamma , \mathcal{T}} $, we put 
\begin{gather*}
\begin{split}
\eta &=  \arccos \left( p\cos \left( \xi + \alpha -\frac{\gamma}{2} \right)  \right) - \arccos \left( p\cos \left( \xi (\theta ) + \alpha -\frac{\gamma}{2} \right)  \right) \\
&= \arccos \left( p\cos \left( \xi + \alpha -\frac{\gamma}{2} \right)  \right) - 2\pi - \frac{\gamma}{2} + \theta  .
\end{split}
\end{gather*}
Note that $\xi ' (\eta ) : = (d\xi / d\eta )(\eta)   $ is smooth and does not vanish in a small neighborhood of $\eta =0$.
Then we have
$$
\lambda_1 (\xi )= e^{i\theta} e^{i\eta} , \quad \theta \in J_{\gamma , 1} \setminus J_{\gamma ,\mathcal{T}} ,
$$
$$
\lambda_2 (\xi )= e^{i\theta} e^{-i\eta} ,\quad \theta \in J_{\gamma ,2} \setminus J_{\gamma ,\mathcal{T}} .
$$
Hence the equation (\ref{S3_eq_helmholtz}) is rewritten in the variable $\eta$ as 
\begin{gather}
\chi (\eta ) P_1 (\eta  ) \widehat{u} (\eta )= \frac{\chi (\eta ) P_1 (\eta ) \widehat{f} (\eta )}{e^{i\theta} e^{i\eta} -e^{iz}} , \quad  \theta \in J_{\gamma , 1} \setminus J_{\gamma ,\mathcal{T}} , \label{S3_eq_helmholtz2} \\
\chi (\eta ) P_2 (\eta  ) \widehat{u} (\eta )= \frac{\chi (\eta ) P_2 (\eta ) \widehat{f} (\eta )}{e^{i\theta} e^{-i\eta} -e^{iz}} , \quad \theta \in J_{\gamma ,2} \setminus J_{\gamma ,\mathcal{T}}  . \label{S3_eq_helmholtz3} 
\end{gather}
Letting $z= \theta -i \log (1\mp \epsilon )$ for small $\epsilon >0$, we consider 
\begin{equation}
\widehat{u}_{1,\pm \epsilon} (\eta )= \frac{ e^{-i\theta} \widehat{f} _{1,\chi} (\eta )}{e^{i\eta} -1 \pm \epsilon } , \quad \widehat{u}_{2,\pm \epsilon} (\eta )= \frac{ e^{-i\theta} \widehat{f} _{2,\chi} (\eta )}{e^{-i\eta} -1 \pm \epsilon } ,
\label{S3_eq_helmholtz4}
\end{equation}
where $\widehat{f} _{j,\chi} = \chi P_j \widehat{f} $.

\begin{lemma}
Let $\theta \in J_{\gamma ,j} \setminus J_{\gamma , \mathcal{T}} $ and $ \widehat{f} \in \mathcal{B} ({\bf T})$.
For any $\widehat{g}\in \mathcal{B} ({\bf T})$, we put $\widehat{g}_{\chi} = \chi \widehat{g} $.
Then there exist the limits $\widehat{u} _{j,\pm} := \lim _{\epsilon \downarrow 0} \widehat{u}_{j,\pm \epsilon} $ in the weak $*$ sense i.e. for any $\widehat{g} \in \mathcal{B} ({\bf T})$
$$
(\widehat{u} _{j,\pm \epsilon} , \widehat{g}_{\chi} ) \to (\widehat{u} _{j,\pm} , \widehat{g}_{ \chi } ) ,
$$
as $\epsilon \downarrow 0$ with the estimates
\begin{equation}
\| \widehat{u} _{j,\pm} \| _{\mathcal{B}^* ({\bf R})} \leq \kappa \| \widehat{f} _{j,\chi} \| _{\mathcal{B} ({\bf R})} ,
\label{S3_eq_BB*}
\end{equation}
for a constant $ \kappa >0$.
Moreover, we have
\begin{gather}
\widetilde{u} _{1,\pm} (s) \mp e^{-i\theta} H(\mp (s+1)) \int_{-\infty}^{\infty} \widetilde{f} _{1,\chi} (t)dt \to 0 , \label{S3_eq_WF1}
\\
\widetilde{u} _{2,\pm} (s) \mp e^{-i\theta} H(\pm (s-1)) \int_{-\infty}^{\infty} \widetilde{f} _{2,\chi} (t)dt \to 0 ,\label{S3_eq_WF2}
\end{gather}
as $|s| \to \infty $, where $H(s)$ is the Heaviside function.

\label{S3_lem_LAP1}
\end{lemma}

Proof.
Let us compute $\widehat{u}_{1,+\epsilon } $.
For other cases, the proofs are parallel.
Since the support of $\widehat{f} _{1,\chi}$ is small, we have
\begin{equation}
(\widehat{u}_{1,+\epsilon} , \widehat{g}_{\chi} )= \frac{e^{-i\theta}}{2\pi} \int _{-\infty}^{\infty} \int _{-\infty}^{\infty} \widetilde{F}_{+,\epsilon} (s-t) \widetilde{f}_{1,\chi}  (t) dt \, \overline{\widetilde{h} _{\chi} (s) }ds ,
\label{S3_eq_LAP11}
\end{equation}
where $\widetilde{F} _{+,\epsilon} $ is defined by (\ref{A_def_Fpme}) and $ \widetilde{h} _{\chi} (s) $ is the Fourier transformation of $ \widehat{g} _{\chi} (\eta ) \xi ' (\eta )$.
Then the equality (\ref{S3_eq_LAP11}) is estimated by 
$$
| (\widehat{u} _{1,+\epsilon} ,\widehat{g} _{\chi} )|\leq \frac{1}{2\pi} \int _{-\infty}^{\infty} |\widetilde{f} _{1, \chi} (t) | dt \int _{-\infty}^{\infty} |\widetilde{h} _{ \chi} (s) | ds ,
$$
in view of Lemmas \ref{A_lem_L2} and \ref{A_lem_fourier}.
The assertion (1) of Lemma \ref{S2_lem_ell1} implies 
$$
| (\widehat{u} _{1,+\epsilon} , \widehat{g} _{\chi} )|\leq \kappa \| \widehat{f} _{1,\chi} \| _{\mathcal{B} ({\bf R})} \| \widehat{g} _{\chi} \| _{\mathcal{B} ({\bf R})} ,
$$
for a constant $\kappa >0$.
This inequality implies $ \| \widehat{u} _{1,+\epsilon} \| _{\mathcal{B}^* ({\bf R})} \leq \kappa \| \widehat{f} _{1,\chi} \| _{\mathcal{B} ({\bf R})} $.
From Lemma \ref{A_lem_fourier}, if $\widehat{f},\widehat{g} \in C^{\infty } ({\bf T} )$, then $(\widehat{u} _{1,+\epsilon} , \widehat{g} _{\chi} ) $ converges to $(\widehat{u} _{1,+} , \widehat{g}_{\chi} )$ as $\epsilon \downarrow 0$, with
\begin{equation}
\widetilde{u} _{1,+} (s)= e^{-i\theta} \left(  \int _{s+1}^{\infty} \widetilde{f}_{1,\chi} (t)dt + (I_+ * \widetilde{f} _{1,\chi} )(s) - (I_- * \widetilde{f}_{1,\chi} )(s)   \right) .
\label{S3_eq_LAP12}
\end{equation}
If $\widehat{f}, \widehat{g}\in \mathcal{B} ({\bf T})$, we approximate them by smooth functions and we can obtain (\ref{S3_eq_BB*}) for $\widehat{u} _{1,+} $.

Note that Lemmas \ref{S2_lem_ell1} and \ref{A_lem_L2} implies $I_{\pm} * \widetilde{f}_{1,\chi} \in L^2 ({\bf R} )$ so that $(I_{\pm} *\widetilde{f}_{1,\chi } )(s)\to 0$ as $|s| \to \infty $.
Then we obtain (\ref{S3_eq_WF1}) for $\widehat{u}_{1,+}$ from (\ref{S3_eq_LAP12}).

Similarly, we can see 
\begin{gather}
\widetilde{u} _{1,-} (s)= e^{-i\theta} \left( - \int_{-\infty}^{s+1} \widetilde{f} _{1,\chi} (t)dt + (I_+ * \widetilde{f}_{1,\chi} )(s) -(I_- * \widetilde{f}_{1,\chi} )(s)     \right) , \label{S3_eq_LAP13} 
\\
\widetilde{u}_{2,+} (s)= e^{-i\theta} \left( \int _{-\infty}^{s-1} \widetilde{f}_{2,\chi } (t)dt +( I^{\dagger}_+ * \widetilde{f}_{2,\chi} )(s) - (I_-^{\dagger} * \widetilde{f}_{2,\chi} )(s)  \right) , \label{S3_eq_LAP14} \\
\widetilde{u}_{2,-} (s)= e^{-i\theta} \left( -\int ^{\infty}_{s-1}\widetilde{f}_{2,\chi } (t)dt +( I^{\dagger}_+ * \widetilde{f}_{2,\chi} )(s) - (I_-^{\dagger} * \widetilde{f}_{2,\chi} )(s)  \right) , \label{S3_eq_LAP15}
\end{gather}
where $I_{\pm}^{\dagger} (s)= I_{\pm} (-s)$.
\qed

\begin{lemma}
For $ \theta \in J_{\gamma ,j} \setminus J_{\gamma , \mathcal{T}} $ and $ \widehat{f} ,\widehat{g} \in \mathcal{B} ({\bf T} )$, we have 
$$
( \widehat{u} _{j,+} -\widehat{u}_{j,-} , \widehat{g} _{j,\chi} )= 2\pi e^{-i\theta}  \widehat{f} _{j,\chi} (0) \cdot \overline{ \widehat{g} _{j,\chi} (0)} \cdot \frac{d\xi }{d\eta } (0) ,
 $$
 where $\widehat{g}_{j,\chi} = \chi P_j \widehat{g} $.
\label{S3_lem_LAP2}
\end{lemma}

Proof.
This lemma is a direct consequence of (\ref{S3_eq_LAP12})-(\ref{S3_eq_LAP15}).
In fact, we have 
$$
( \widehat{u} _{j,+} -\widehat{u}_{j,-} , \widehat{g} _{j,\chi} )= e^{-i\theta}  \int_{-\infty}^{\infty} \widetilde{f}_{j,\chi} (s) ds \overline{ \int_{-\infty}^{\infty} \widetilde{h}_{j,\chi} (s) ds} ,
$$
where $\widetilde{h}_{\chi} (s)$ has been used in (\ref{S3_eq_LAP11}).
This equality shows the lemma.
\qed

\medskip

Next we consider a radiation condition which guarantees the uniqueness of the solution to the equation
\begin{gather}
( \lambda_1 (\xi )-e^{i\theta} )\chi (\xi ) P_1 (\xi ) \widehat{u} (\xi )=\chi (\xi ) P_1 (\xi )  \widehat{f} (\xi ) , \label{S3_eq_helmholtz_rad1} \\
 ( \lambda_2 (\xi )-e^{i\theta} )\chi (\xi ) P_2 (\xi ) \widehat{u} (\xi )= \chi (\xi ) P_2 (\xi ) \widehat{f} (\xi ),
\label{S3_eq_helmholtz_rad2}
\end{gather}
for $\widehat{f}\in \mathcal{B} ({\bf T} )$ and $ \theta \in J_{\gamma} \setminus J_{\gamma , \mathcal{T}} $.

\begin{definition}
(1) For $ u\in \mathcal{S} ' ({\bf R}) $, the wave front set $WF^* (u)$ is defined as follows.
For $  (\eta_0 , \omega ) \in {\bf R} \times \{ \pm 1 \} $, $(\eta_0 , \omega ) \not\in WF^* (u)$ if there exists $\psi \in C_0^{\infty} ({\bf R})$ such that $\psi (\eta_0 )=1$ and
\begin{equation}
\lim _{R\to \infty} \frac{1}{R} \int _{|s|<R} | C_{\omega} (s) \widetilde{\psi u} (s) |^2 ds =0 ,
\label{S3_defeq_WF}
\end{equation}
where $ C_{\omega} (s)$ is the characteristic function of the set $\{ s\in {\bf R} \ ; \ \mathrm{sgn} (s)=\omega \} $. \\
(2) Let $u$ be a distribution on ${\bf T}$.
For $(\eta_0 , \omega )\in {\bf T} \times \{ \pm 1 \} $, $(\eta_0 ,\omega ) \not\in WF^* (u)$ if there exists $\psi \in C^{\infty} ({\bf T})$ such that $\psi (\eta_0 )=1$ with small support and $\widetilde{\psi u}$ satisfies (\ref{S3_defeq_WF}).
\label{S3_def_WF}
\end{definition}

Now let us show the uniqueness of solutions $\widehat{u} _{\pm} \in \mathcal{B}^* ({\bf T} )$ to the equation (\ref{S3_eq_helmholtz_rad1})-(\ref{S3_eq_helmholtz_rad2}) with the radiation condition $WF^* (\chi P_1 \widehat{u} _{\pm} ) = \{ (0,\mp 1) \} $  for $\theta \in J_{\gamma ,1} \setminus J_{\gamma , \mathcal{T}} $ or $WF^* ( \chi P_2 \widehat{u} _{\pm })=\{ (0,\pm 1)\} $ for $\theta \in J_{\gamma ,2} \setminus J_{\gamma , \mathcal{T}} $.

\begin{lemma}
(1) For $ \widehat{u}_{j,\pm} $ given in Lemma \ref{S3_lem_LAP1}, we have 
$$
WF^* (\widehat{u} _{1,\pm} )=\{ (0, \mp 1)\} ,\quad WF^* (\widehat{u} _{2,\pm} )=\{ (0, \pm 1)\} .
$$
Moreover, $\widehat{u} _{1,\pm} $ and $\widehat{u} _{2,\pm} $ satisfy
\begin{gather}
\widehat{u} _{1,\pm} - \frac{ e^{-i\theta} e^{-i\eta} }{i\eta \pm 0} \otimes \widehat{f} _{1,\chi} (0) \in \mathcal{B}_0^* ({\bf R}) , \label{S3_eq_radcondition11} \\
\widehat{u} _{2,\pm} + \frac{ e^{-i\theta} e^{i\eta} }{i\eta \mp 0} \otimes \widehat{f} _{2,\chi} (0) \in \mathcal{B}_0^* ({\bf R}) . \label{S3_eq_radcondition12} 
\end{gather}
(2) Let $\widehat{u}\in \mathcal{B}^* ({\bf T})$ be a solution of (\ref{S3_eq_helmholtz_rad1})-(\ref{S3_eq_helmholtz_rad2}) for $\theta \in J_{\gamma ,1} \setminus J_{\gamma , \mathcal{T}} $. 
The solution $\widehat{u}$ satisfies $WF^* ( \chi P_1 \widehat{u} )= \{ (0, \mp 1)\} $ if and only if $\chi P_1 \widehat{u}=\widehat{u} _{1,\pm} $. 
When $\theta \in J_{\gamma ,2} \setminus J_{\gamma , \mathcal{T}} $, the solution $\widehat{u}$ satisfies $WF^* (\chi P_2 \widehat{u} )=\{ (0, \pm 1)\} $ if and only if $\chi P_2  \widehat{u} = \widehat{u} _{2,\pm} $.
\label{S3_lem_radcondition}
\end{lemma}

Proof.
Note that (\ref{S3_eq_radcondition11}) and (\ref{S3_eq_radcondition12}) are direct consequences of (\ref{S3_eq_WF1}) and (\ref{S3_eq_WF2}).
For the proof of the wave front sets, we have only to prove for $\widehat{f}\in C^{\infty} ({\bf T})$.
Since we have 
$$
\widetilde{u} _{1,+} (s)= e^{-i\theta} \int _{s+1}^{\infty} \widetilde{f} _{1,\chi } (t)dt , \quad \widetilde{u}_{1,-} (s)=- e^{-i\theta} \int _{-\infty}^{s+1} \widetilde{f}_{1,\chi} (t)dt ,
$$
up to terms in $L^2 ({\bf R})$, we see 
$$
|\widetilde{u} _{1,+} (s) |\leq \kappa _n \int _{s+1}^{\infty} (1+|t|)^{-n} dt , \quad |\widetilde{u} _{1,-} (s) |\leq \kappa _n \int ^{s+1}_{-\infty} (1+|t|)^{-n} dt ,
$$
for any $n>0$ with constants $ \kappa _n >0$.
Therefore, we can see $WF^* (\widehat{u} _{1,\pm} )=\{ (0,\mp 1 )\} $.
Noting that
$$
\widetilde{u} _{2,+} (s)= e^{-i\theta} \int_{-\infty}^{s-1} \widetilde{f}_{2,\chi} (t)dt ,\quad \widetilde{u}_{2,-} (s)=-e^{-i\theta} \int_{s-1}^{\infty} \widetilde{f}_{2,\chi} (t)dt ,
$$
up to terms in $L^2 ({\bf R})$, we also have $WF^* (\widehat{u} _{2,\pm} )=\{ (0, \pm 1 ) \} $.

Let us turn to the proof of the assertion (2).
For $\theta \in J_{\gamma ,1} \setminus J_{\gamma , \mathcal{T} }$, let $ \widehat{u} \in \mathcal{B}^* ({\bf T})$ satisfy the equation (\ref{S3_eq_helmholtz_rad1})-(\ref{S3_eq_helmholtz_rad2}) and the condition $WF^* (\chi P_1 \widehat{u})= \{ (0,-1) \} $.
In view of Lemma \ref{S3_lem_distinct}, we have only to consider (\ref{S3_eq_helmholtz_rad1}).
Then $\widehat{v} := \widehat{u}_{1,+} -\chi P_1 \widehat{u} $ satisfies the equation
\begin{equation}
e^{i\theta} (e^{i\eta} -1 ) \widehat{v} (\eta )=0 ,
\label{S3_eq_homogeneous}
\end{equation}
in the variable $\eta $ with the condition $WF^* (\chi P_1 \widehat{v} )=\{ (0,-1) \} $.
Passing to the Fourier transformation, the equation (\ref{S3_eq_homogeneous}) implies that $\widetilde{v} (s)=\widetilde{v} (s-1)$ for any $s\in {\bf R}$ i.e. $\widetilde{v} $ is 1-periodic.
Thus it follows from $WF^* (\widehat{v} )=\{ (0,-1) \} $ that
$$
\int _{0<s<1} |\widetilde{v} (s)|^2 ds = \frac{1}{n} \int _{0<s<n} |\widetilde{v} (s)|^2 ds , 
$$
for positive integers $n$ converges to zero as $n\to\infty $.
We obtain $\widetilde{v} =0$ i.e. $\chi P_1 \widehat{u} = \widehat{u} _{1,+} $.

For $\theta \in J_{\gamma ,2} \setminus J_{\gamma , \mathcal{T}} $, we can see that the uniqueness of the solution of the equation (\ref{S3_eq_helmholtz_rad1})-(\ref{S3_eq_helmholtz_rad2}) with the condition $WF^* (\chi P_2 \widehat{u} )= \{ (0,1) \} $ by the same way.
\qed

\medskip

Summarizing the argument in this subsection, we have proven the limiting absorption of $\widehat{R}_0 (z)$ as follows.

\begin{theorem}
(1) For $\theta \in J_{\gamma } \setminus J_{\gamma , \mathcal{T}} $ and $\widehat{f}\in \mathcal{B} ({\bf T})$, there exists a weak $*$ limit $\widehat{R}_0 (\theta \pm i0) \widehat{f} := \lim_{\epsilon \downarrow 0} \widehat{R}_0 (\theta -i\log (1\mp \epsilon )) \widehat{f} $ in the sense
$$
( \widehat{R}_0 (\theta -i\log (1\mp \epsilon )) \widehat{f} , \widehat{g} ) \to (\widehat{R}_0 (\theta \pm i0) \widehat{f} ,\widehat{g} ) ,  
$$
as $\epsilon \downarrow 0$ for any $\widehat{g} \in \mathcal{B} ({\bf T} )$. \\
(2) $ \widehat{R}_0 (\theta \pm i0) \widehat{f} $ for $\widehat{f} \in \mathcal{B} ({\bf T})$ satisfies 
$$
\| \widehat{R}_0 (\theta \pm i0 )\widehat{f} \| _{\mathcal{B}^* ({\bf T})} \leq \kappa \| \widehat{f} \| _{\mathcal{B} ({\bf T} )} ,
$$ 
where the constant $\kappa >0$ is independent of $\theta$ if $\theta $ varies over a compact interval in $ J_{\gamma} \setminus J_{\gamma , \mathcal{T} } $. \\
(3) The mapping $ J_{\gamma} \setminus J_{\gamma , \mathcal{T}} \ni \theta \mapsto (\widehat{R}_0 (\theta \pm i0) \widehat{f} ,\widehat{g} ) $ is continuous  for $\widehat{f} , \widehat{g} \in \mathcal{B} ({\bf T})$. \\
(4) $\widehat{u} _{\pm} = \widehat{R}_0 (\theta \pm i0)\widehat{f}$ for $\widehat{f}\in \mathcal{B} ({\bf T})$ are the unique solutions of $ ( \widehat{U}_0 -e^{i\theta} ) \widehat{u}_{\pm} = \widehat{f} $ satisfying the condition
$$
WF^* ( P_1 \widehat{u} _{\pm} ) = \{ (\xi (\theta ) , \mp 1) \} , \quad WF^* (P_2 \widehat{u}_{\pm} )=\{ (\xi (\theta ) ,\pm 1 )\} ,
$$
for $\xi (\theta ) \in M( \theta )$, and 
\begin{gather*}
P_1 \widehat{u} _{\pm} = \frac{e^{-i \gamma /2 } \mathrm{exp} (-i \arccos (p\cos (\xi + \alpha -\gamma /2))  )}{i \left( \arccos (p\cos (\xi + \alpha -\gamma /2)) - \theta +\gamma /2 \right) \pm 0 } \otimes ( P_1 \widehat{f} )(\xi (\theta ))   , \\
P_2 \widehat{u} _{\pm} =- \frac{e^{-i\gamma /2} \mathrm{exp} (i\arccos (p\cos (\xi +\alpha -\gamma /2)))}{i \left( \arccos (p\cos (\xi +\alpha -\gamma /2 )) -2\pi + \theta -\gamma /2 \right) \mp 0 } \otimes  ( P_2 \widehat{f} )(\xi (\theta ))   ,
\end{gather*}
up to terms in $ \mathcal{B}_0^* ({\bf T} )$. \\
(5) For $ \theta \in J_{\gamma ,j} \setminus J_{\gamma , \mathcal{T} } $ and $\widehat{f},\widehat{g} \in \mathcal{B} ({\bf T})$, we have 
\begin{gather}
\begin{split}
&( ( \widehat{R} _0 (\theta +i0) - \widehat{R} _0 (\theta -i0))  \widehat{f} ,  \widehat{g} ) \\
&= 2\pi e^{-i\theta} \sum _{\xi (\theta ) \in M(\theta )}  ( (P_j \widehat{f} )(\xi (\theta )) ,(P_j \widehat{g} )(\xi (\theta )) )_{{\bf C}^2}  \frac{d\xi}{d\theta} (\theta )  , 
\end{split}
\label{S3_eq_parsevalformulaxi}
\end{gather}
for $\xi (\theta )\in M(\theta )$. 
\label{S3_thm_LAPR0}
\end{theorem}

Proof.
We have only to make a remark for the assertion (5).
Let $ \widehat{\phi} (\xi )$ be a suitable function such that its support is in a small neighborhood of a point $\xi (\theta ) \in M(\theta )$.
The integration of $\widehat{\phi} (\xi )$ is rewritten in the variable $\eta$ as  
$$
\int _{{\bf T} } \widehat{\phi} (\xi )d\xi = \int _{-\delta }^{\delta }  \widehat{\phi} (\eta )  \xi ' (\eta)  d \eta ,
$$
for a small $\delta >0$.
By the definition of the variable $\eta$, we have 
$$
\frac{d\xi }{d\eta } ( \eta )= \frac{\sin (\eta + K_j ) }{\sqrt{ p^2 - \cos^2 (\eta +K_j )} } , 
$$
where $K_1 = \theta - \gamma /2$ for $\theta \in J_{\gamma ,1} \setminus J_{\gamma , \mathcal{T}} $ and $K_2 = -\theta + \gamma /2 $ for $ \theta \in J_{\gamma ,2} \setminus J_{\gamma , \mathcal{T}} $.
In view of (\ref{S2_eq_determinant}), we also have 
$$
\frac{d\xi }{d\theta } (\theta )= \frac{\sin K_j }{\sqrt{p^2 - \cos^2  K_j }} >0 ,
$$
for $\theta \in J_{\gamma , j} \setminus J_{\gamma , \mathcal{T} } $.
Thus we obtain
$$
\frac{d\xi }{d\eta } (0) = \frac{d\xi}{d\theta } (\theta ).
$$
Therefore, Lemma \ref{S3_lem_LAP2} can be rewritten in the variable $\xi$ as (\ref{S3_eq_parsevalformulaxi}).
\qed

\medskip

Let us give some remarks for adjoint operators $ R_0 (\theta \pm i0 )^* $. 
By the similar way, we can show that the limits $ R_0 (\theta \pm i0 )^* := \lim _{\epsilon \downarrow 0} R_0 (\theta - i \log  ( 1\mp \epsilon )  )^* \in {\bf B} ( \mathcal{B} ({\bf Z}) ; \mathcal{B}^* ({\bf Z})) $ exist for $\theta \in J_{\gamma} \setminus J_{\gamma , \mathcal{T} }$ in the sense
$$
( f, R_0 (\theta - i \log  ( 1\mp \epsilon )  )^* g) \to ( f, R_0 (\theta \pm i 0 ) ^* g)  , \quad f,g \in \mathcal{B} ({\bf Z} ) .
$$
By the definition, we have
$$
( R_0 ( \theta \pm i0 ) f,g)= (f, R_0 (\theta \pm i0 )^* g) ,\quad f,g \in \mathcal{B} ({\bf Z} ).
$$
The assertions (1)-(4) in Theorem \ref{S3_thm_LAPlattice} hold for the adjoint operators $ R_0 (\theta \pm i0)^* $, noting that 
\begin{equation}
WF^* ( P_1 \mathcal{U} u_{\pm}^* )= \{ (\xi (\theta ) , \pm 1 )\} , \quad WF^* ( P_2 \mathcal{U} u _{\pm} ^* )=\{ (\xi (\theta ) , \mp 1 )\} ,
\label{S3_eq_WFadjoint}
\end{equation}
where $u _{\pm}^* = R_0 (\theta \pm i0 )^* f$ for $f\in \mathcal{B} ({\bf Z} )$.


\subsection{Spectral representation}

We define a generalized Fourier transformation associated with $U_0$.
First we introduce a Hilbert space on ${\bf T}$.
We can see that
\begin{gather}
\begin{split}
\int_{{\bf T}} |\widehat{f} (\xi )|^2 _{{\bf C}^2} d\xi &= \sum_{j=1}^2 \int_{{\bf T}} | P_j (\xi ) \widehat{f} (\xi )|^2_{{\bf C}^2} d\xi \\
&= \sum _{j=1}^2 \int _{J_{\gamma ,j}} \sum _{\xi (\theta )\in M(\theta )} | P_j (\xi (\theta )) \widehat{f} (\xi (\theta )) |^2_{{\bf C}^2}  \frac{d\xi}{d\theta} (\theta ) d\theta .
\end{split}
\label{S3_eq_L2toH}
\end{gather}
This observation motivates us to introduce the following Hilbert spaces.
Let $\widetilde{{\bf h}} (\theta )$ is the space of ${\bf C}$-valued functions on $M(\theta )$ with its inner products
$$
(\phi , \psi )_{\widetilde{{\bf h}} (\theta )} = \sum_{\xi (\theta )\in M(\theta )}  \phi (\xi (\theta )) \cdot \overline{ \psi (\xi (\theta )) } \cdot  \frac{d\xi}{d\theta} (\theta ) ,
$$
for $ \phi , \psi : M(\theta ) \to {\bf C} $.
Let $a_j (\xi )$ be a normalized eigenvector of $\widehat{U}_0 (\xi )$ with eigenvalue $\lambda_j (\xi )$.
Since $\widehat{U}_0 (\xi ) $ is unitary, the projection $P_j (\xi )$ is given by $ P_j (\xi )\widehat{f} (\xi )= (\widehat{f} (\xi ) , a_j (\xi ))_{{\bf C}^2} a_j (\xi )$.
Then we define
$$
{\bf H}_j = L^2 (J_{\gamma ,j} ;  \widetilde{{\bf h}} (\theta ) a_j ; d\theta ) ,
$$
where $ \widetilde{{\bf h}} (\theta ) a_j =\{  \phi a_j |_{M(\theta )} \ ; \ \phi \in \widetilde{{\bf h}} (\theta )\} $.
Note that the eigenfunction $a_j$ cannot chosen uniquely.
In the strict sense, ${\bf H}_j$ is the equivalence class  in view of the equivalent relation 
$$
( \widehat{f} , a_j ) \sim ( \widehat{g} , b_j ) \iff (\widehat{f} (\xi) ,a_j (\xi ) )_{{\bf C}^2} a_j (\xi )=(\widehat{g} (\xi) ,b_j (\xi ) )_{{\bf C}^2} b_j (\xi ) ,
$$
where $\widehat{f},\widehat{g} \in L^2 ({\bf T} )$ and $a_j (\xi ) ,b_j (\xi)$ are normalized eigenvectors of $\widehat{U}_0 (\xi )$ with eigenvalue $\lambda_j (\xi )$.
In the following, we put 
$$
{\bf H} = \oplus _{j=1}^2 {\bf H}_j ,\quad {\bf h} (\theta )= \oplus _{j=1}^2  \widetilde{{\bf h}} (\theta ) a_j .
$$
Note that, for $ \phi , \psi \in {\bf h} (\theta )$, there exist $ \phi_j , \psi _j \in \widetilde{{\bf h}}  (\theta )$ such that $ \phi = \oplus _{j=1}^2 \phi_j a_j |_{M(\theta )}$ and $ \psi =\oplus _{j=1}^2 \psi_j a_j |_{M(\theta )} $.
Moreover, we have
$$
(\phi , \psi )_{{\bf h} (\theta )} = \sum_{j=1}^2 \sum _{\xi (\theta ) \in M(\theta )} \phi_j (\xi (\theta )) \cdot \overline{\psi_j (\xi (\theta ))} \cdot \frac{d\xi}{d\theta} (\theta ) .
$$

We define the operator $\widehat{\mathcal{F}}_0 (\theta ) =(\widehat{\mathcal{F}}_{0,1} (\theta ),\widehat{\mathcal{F}}_{0,2} (\theta ))$ by
\begin{equation}
\widehat{\mathcal{F}}_{0,j} (\theta ) \widehat{f} =\left\{
\begin{split}
P_j  \widehat{f} |_{M(\theta )}  &, \quad \theta \in J_{\gamma ,j} ,\\
0 &, \quad \theta \not\in J_{\gamma ,j} ,
\end{split}
\right. \quad \widehat{f} \in \mathcal{B} ({\bf T} ).
\label{S3_def_F0hat}
\end{equation}
Let $(\widehat{\mathcal{F}}_0  \widehat{f} )(\theta )= \widehat{\mathcal{F}}_0 (\theta )\widehat{f} $.
Then we have the following fact.

\begin{lemma}
The operator $\widehat{\mathcal{F}}_0$ can be extended uniquely to a unitary operator from $L^2 ({\bf T} ; {\bf C}^2 ) $ to ${\bf H} $. 
Moreover, we have $( \widehat{\mathcal{F}}_0 \widehat{U} _0 \widehat{f} )(\theta ) = e^{i\theta} (\widehat{\mathcal{F}}_0 \widehat{f})(\theta )$ and $( \widehat{\mathcal{F}}_0 \widehat{U} _0^* \widehat{f} )(\theta ) = e^{-i\theta} (\widehat{\mathcal{F}}_0 \widehat{f})(\theta )$ for any $\widehat{f} \in L^2 ({\bf T}  ; {\bf C}^2 ) $.
\label{S3_lem_fourier0}
\end{lemma}

Proof. 
We have only to verify $ \widehat{\mathcal{F}}_0 (\theta ) \widehat{U}_0^* \widehat{f} = e^{-i\theta} \widehat{\mathcal{F}} _0 (\theta )\widehat{f} $.
In fact, we can see 
\begin{gather*}
\begin{split}
p^* (\xi , \theta )  := &\, \mathrm{det} (\widehat{U}_0 (\xi )^* -e^{-i\theta })\\
 = & \, 2p e^{-i(\theta + \gamma/2 )} \left( - \cos \left( \xi + \alpha - \frac{\gamma}{2} \right) + \frac{1}{p} \cos \left( -\theta +\frac{\gamma}{2} \right)   \right) .
\end{split}
\end{gather*}
If $\xi (\theta )\in M(\theta )$, it follows $p^* (\xi (\theta ), \theta ) =0$ from $ p \cos ( \xi (\theta )+ \alpha -\gamma /2)= \cos (\theta - \gamma /2)$. 
Since $ \widehat{U}_0 (\xi )\in U(2)$ for every $\xi \in {\bf T} $, we have 
$$
\widehat{U}_0 (\xi ) P_j (\xi ) \widehat{f} (\xi )= e^{i\theta} P_j (\xi ) \widehat{f} (\xi ), 
$$
if and only if 
$$
\widehat{U}_0 (\xi )^* P_j (\xi ) \widehat{f} (\xi )= e^{-i\theta} P_j (\xi ) \widehat{f} (\xi ).
$$
Thus we obtain the lemma.
\qed

\medskip

The Parseval formula for $\widehat{\mathcal{F}}_0 $ is given as follows.

\begin{lemma}
For $\theta \in J_{\gamma } \setminus J_{\gamma , \mathcal{T}} $, we have 
$$
( \widehat{R} _0 (\theta +i0) \widehat{f} - \widehat{R} _0 (\theta -i0)  \widehat{f} ,  \widehat{g} )= 2\pi e^{-i\theta} ( \widehat{\mathcal{F}}_{0} (\theta ) \widehat{f}    , \widehat{\mathcal{F}}_{0} (\theta ) \widehat{g}  )_{{\bf h} (\theta )} ,
$$
for $\widehat{f} , \widehat{g} \in \mathcal{B} ({\bf T})$.
It follows from the above equality that $\widehat{\mathcal{F}}_0 (\theta )\in {\bf B} ( \mathcal{B} ({\bf T}) ;  {\bf h} (\theta ))$  with the estimate
$$
\| \widehat{\mathcal{F}} _{0} (\theta )\widehat{f} \|_{{\bf h}_j (\theta )} \leq \kappa \| \widehat{f} \| _{\mathcal{B} ({\bf T})} , \quad \widehat{f} \in \mathcal{B} ({\bf T}) ,
$$
for a constant $\kappa >0$.
\label{S3_lem_parseval}
\end{lemma}

Proof.
This lemma is also a direct consequence of (\ref{S3_eq_LAP12})-(\ref{S3_eq_LAP15}) and Lemma \ref{S3_thm_LAPR0}.
\qed

\medskip

Then we have $\widehat{\mathcal{F}}_0 (\theta )^* \in {\bf B} ({\bf h} (\theta ) ;\mathcal{B}^* ({\bf T} ))$. 
The equality $ \widehat{\mathcal{F}}_0 (\theta ) ( \widehat{U}_0^* - e^{-i \theta } ) \widehat{f} =0$ for $\widehat{f}\in \mathcal{B} ({\bf T})$ implies 
$$
(\widehat{U}_0 -e^{i\theta } ) \widehat{\mathcal{F}}_0 (\theta )^* \phi =0 , \quad \phi \in {\bf h} (\theta ) .
$$
Therefore, $\widehat{\mathcal{F}} _0 (\theta )^* $ is an eigenoperator of $\widehat{U}_0 $.


\subsection{On Lattice}

Passing to the inverse Fourier transformation $\mathcal{U}^* $, the arguments in \S 3.1-3.2 are translated on the lattice ${\bf Z}$. 
The following theorem is a direct consequence. 
We omit their proofs.

In the following, we define
$$
\mathcal{F}_0 = \widehat{\mathcal{F}}_0 \mathcal{U}  , \quad \mathcal{F}_0 (\theta ) = \widehat{\mathcal{F}}_0 (\theta ) \mathcal{U} , \quad \mathcal{F}_{0,j} (\theta ) = \widehat{\mathcal{F}}_{0,j} (\theta ) \mathcal{U} .
$$

\begin{theorem}
(1) For $\theta \in J_{\gamma } \setminus J_{\gamma , \mathcal{T}} $ and $f\in \mathcal{B} ({\bf Z})$, there exists a weak $*$ limit $R_0 (\theta \pm i0) f := \lim_{\epsilon \downarrow 0} R_0 (\theta -i\log (1\mp \epsilon )) f $ in the sense
$$
( R_0 (\theta -i\log (1\mp \epsilon )) f , g ) \to (R_0 (\theta \pm i0) f ,g ) ,  
$$
as $\epsilon \downarrow 0$ for any $g \in \mathcal{B} ({\bf Z} )$. \\
(2) $ R_0 (\theta \pm i0) f $ for $f \in \mathcal{B} ({\bf Z})$ satisfies 
$$
\| R_0 (\theta \pm i0 )f \| _{\mathcal{B}^* ({\bf Z})} \leq \kappa \| f \| _{\mathcal{B} ({\bf Z} )} ,
$$ 
where the constant $\kappa >0$ is independent of $\theta$ if $\theta $ varies over a compact interval in $ J_{\gamma} \setminus J_{\gamma , \mathcal{T} } $. \\
(3) The mapping $ J_{\gamma} \setminus J_{\gamma , \mathcal{T}} \ni \theta \mapsto (R_0 (\theta \pm i0) f ,g ) $ is continuous  for $f , g \in \mathcal{B} ({\bf Z})$. \\
(4) $u _{\pm} = R_0 (\theta \pm i0)f$ for $f\in \mathcal{B} ({\bf Z})$ are the unique solutions of $ (U_0 -e^{i\theta} )u_{\pm} =f$ satisfying the condition  
\begin{equation}
WF^* ( P_1 \mathcal{U} u_{\pm} ) = \{ (\xi (\theta ) , \mp 1) \} , \quad WF^* (P_2 \mathcal{U} u_{\pm} )=\{ (\xi (\theta ) ,\pm 1 )\} ,
\label{S3_eq_radconditionlattice}
\end{equation}
for $\xi (\theta ) \in M( \theta )$, and 
\begin{gather*}
P_1 \mathcal{U} u _{\pm} = \frac{e^{-i \gamma /2 } \mathrm{exp} (-i \arccos (p\cos (\xi + \alpha -\gamma /2))  )}{i \left( \arccos (p\cos (\xi + \alpha -\gamma /2)) - \theta +\gamma /2 \right) \pm 0 } \otimes ( P_1 \mathcal{U} f )(\xi (\theta ))   , \\
P_2 \mathcal{U} u _{\pm} =- \frac{e^{-i\gamma /2} \mathrm{exp} (i\arccos (p\cos (\xi +\alpha -\gamma /2)))}{i \left( \arccos (p\cos (\xi +\alpha -\gamma /2 )) -2\pi + \theta -\gamma /2 \right) \mp 0 } \otimes  ( P_2 \mathcal{U} f)(\xi (\theta ))   ,
\end{gather*}
up to terms in $ \mathcal{B}_0^* ({\bf T} )$. \\
(5) For $ \theta \in J_{\gamma } \setminus J_{\gamma , \mathcal{T} } $, we have 
$$
( R _0 (\theta +i0)  f - R_0 (\theta -i0) f,  g )= 2\pi e^{-i\theta} (\mathcal{F}_{0} (\theta )f ,\mathcal{F}_{0} (\theta )g )_{{\bf h}(\theta )} ,
$$
for $f,g\in \mathcal{B} ({\bf Z})$.
\\
(6) We have $\mathcal{F}_0 (\theta ) \in {\bf B} (\mathcal{B} ({\bf Z}) ; {\bf h} (\theta ))$ and $ \mathcal{F}_0 (\theta )^* \in {\bf B} ({\bf h} (\theta ); \mathcal{B}^* ({\bf Z}))$.
$ \mathcal{F}_0 (\theta )^* \phi $ for $\theta \in J_{\gamma} \setminus J_{\gamma , \mathcal{T}} $ and $\phi \in {\bf h} (\theta )$ satisfies the equations $ (U_0 -e^{i\theta} )\mathcal{F}_0 (\theta )^* \phi =0$ and $(U_0^* -e^{-i\theta} ) \mathcal{F}_0 (\theta )^* \phi =0 $.
\label{S3_thm_LAPlattice}
\end{theorem}

\subsection{Generalized eigenfunction for position-independent QW}

Let $ \delta (\xi - \xi (\theta ))$ be the distribution defined by
$$
\int _{{\bf T}} \delta (\xi - \xi (\theta )) f(x) dx = f(\xi (\theta )) , \quad f\in C({\bf T}) .
$$
In view of $ (\widehat{\mathcal{F}}_0 (\theta )\widehat{f} , \phi )_{{\bf h} (\theta )}= (\widehat{f} , \widehat{\mathcal{F}} _0 (\theta )^* \phi ) $ for $\widehat{f} \in  C ({\bf T}) $ and $\phi \in {\bf h} (\theta )$, the operator $\widehat{\mathcal{F}}_0 (\theta )^* $ defines the distribution for $ \theta \in J_{\gamma ,j} \setminus J_{\gamma , \mathcal{T}} $
\begin{equation}
(\widehat{\mathcal{F}}_0 (\theta )^* \phi )(\xi ) = \sum_{\xi (\theta )\in M(\theta )}  \frac{d\xi }{d\theta } (\theta )   \phi_j (\xi (\theta )) \delta (\xi -\xi (\theta ))  a_j (\xi (\theta ))  ,
\label{S3_def_f0*} 
\end{equation}
where $ \phi = \oplus _{j=1}^2 \phi_j a_j $ for a normalized eigenvector $a_j (\xi )$ of $\widehat{U}_0 (\xi )$.
If $f\in \mathcal{B} ({\bf Z})$, Lemma \ref{S2_lem_ell1} implies $f\in \ell^1 ({\bf Z} )$. 
Then the Fourier series $\widehat{f} (\xi )= (2\pi )^{-1/2} \sum _{x\in {\bf Z}} e^{-ix\xi} f (x)$ converges uniformly so that $\widehat{f} \in C ({\bf T})$.
Thus the generalized eigenfunction of $ U_0$ is given by
\begin{equation}
(\mathcal{F}_0 (\theta )^* \phi )(x) = \frac{1}{\sqrt{2\pi}} \sum_{\xi (\theta )\in M(\theta )}  \frac{d\xi }{d\theta } (\theta )  e^{ix\xi (\theta )}   \phi_j (\xi (\theta ))   a_j (\xi (\theta )) ,
\label{S3_eq_generalizedef0}
\end{equation}
for $ \theta \in J_{\gamma ,j} \setminus J_{\gamma , \mathcal{T}}$.

\begin{lemma}
For $ \theta \in J_{\gamma ,j} \setminus J_{\gamma , \mathcal{T}} $, we have 
$$
\limsup_{R\to \infty} \frac{1}{R} \sum _{|x|<R} | (\mathcal{F}_0 (\theta )^* \phi )(x)|^2 _{{\bf C}^2} \geq \kappa (\theta ) \| \phi \|^2 _{{\bf h} (\theta )} ,
$$
where $ \kappa (\theta )>0$ is a constant.
\label{S3_lem_F0*limit}
\end{lemma}

Proof.
Let $ M(\theta ) = \{ \xi _1 , \xi_2 \} $. 
In view of $\theta \in J_{\gamma ,j} \setminus J_{\gamma , \mathcal{T}} $, we can assume $\xi_1 - \xi_2 \not= 0$ modulo $2\pi $.
Thus we have
\begin{gather}
\begin{split}
| (\mathcal{F}_0 (\theta )^* \phi )(x)|^2 _{{\bf C}^2} = & \,  \theta '(\xi_1 ) ^{-2} |\phi_j (\xi_1 )|^2 + \theta '(\xi_2 ) ^{-2} |\phi_j (\xi_2 )|^2  \\
&+ 2\mathrm{Re} \left( \phi_j (\xi_1 ) \overline{\phi_j (\xi_2 )} e^{ix(\xi_1 -\xi_2 )} \right) \theta '(\xi_1 )^{-1} \theta '(\xi_2 )^{-1}   ,
\end{split}
\label{S3_eq_lim11}
\end{gather}
where $\theta ' (\xi_k )= ( d\theta / d\xi ) | _{\xi= \xi_k} $.
If $\xi _1 - \xi_2 \not= \pi $ modulo $2\pi $, we have 
$$
\sum _{|x|<R} e^{ix(\xi_1 - \xi_2)} = \frac{1-e^{iR(\xi_1 -\xi_2 )} }{1-e^{i (\xi_1 -\xi_2 )}} + \frac{ e^{-i(\xi_1 -\xi_2 )} (1-e^{-i(R-1) (\xi_1 - \xi_2 )})}{1- e^{-i(\xi_1 -\xi_2 )} } ,
$$
for large positive integers $R>1$.
Then $\sum_{|x|<R} e^{ix(\xi_1 - \xi_2)}$ is bounded with respect to $R$ so that we have 
\begin{equation}
\lim _{R\to \infty} \frac{1}{R} \sum_{|x|<R} e^{ix(\xi_1 - \xi_2)} =0 .
\label{S3_eq_lim0}
\end{equation}
If $\xi _1 - \xi_2 =\pi $ modulo $2\pi$, we have 
$$
\sum _{|x|<R} e^{ix(\xi_1 - \xi_2)} = \sum _{|x|<R} (-1)^x =1+ 2\sum _{x=1}^{R-1} (-1)^x .
$$
In this case, $\sum_{|x|<R} e^{ix(\xi_1 - \xi_2)}$ is also bounded with respect to $R$.
Then (\ref{S3_eq_lim0}) holds for $\xi_1 - \xi_2 =\pi $ modulo $2\pi $.
Plugging (\ref{S3_eq_lim11}) and (\ref{S3_eq_lim0}), we have
$$
\limsup_{R\to \infty} \frac{1}{R} \sum_{|x|<R} | (\mathcal{F}_0 (\theta )^* \phi )(x)|^2 _{{\bf C}^2} = \frac{1}{\pi} \left(  \theta '(\xi_1 ) ^{-2} |\phi_j (\xi_1 )|^2 + \theta '(\xi_2 ) ^{-2} |\phi_j (\xi_2 )|^2 \right) .
$$
Taking $ \kappa (\theta )=\pi^{-1} \min \{ \theta '(\xi_1 )^{-1} , \theta ' (\xi_2 )^{-1} \} $, we obtain the lemma.
\qed

\medskip

We can characterize the set of generalized eigenfunctions of $U_0$ by using $ \mathcal{F}_0 (\theta )^* $.

\begin{theorem}
Let $ \theta \in J_{\gamma} \setminus J_{\gamma , \mathcal{T}} $.
Then we have $\mathcal{F}_0 (\theta ) \mathcal{B} ({\bf Z} )= {\bf h} (\theta )$ and $\{ u\in \mathcal{B}^* ({\bf Z} ) \ ; \ (U_0 -e^{i\theta} )u=0 \} = \mathcal{F}_0 (\theta )^* {\bf h} (\theta )$.
\label{S3_thm_rangeF0*}
\end{theorem}

Proof.
In view of Lemmas \ref{S3_lem_parseval} and \ref{S3_lem_F0*limit} and the inequality
\begin{equation}
\limsup _{R\to \infty} \frac{1}{R} \sum _{|x|<R} | u (x)|^2_{{\bf C}^2} \leq \sup _{R>1} \frac{1}{R} \sum _{|x|<R} |u(x)|^2_{{\bf C}^2} ,  \quad u\in \mathcal{B} ({\bf Z}),
\label{S3_eq_ineqB*limsup}
\end{equation}
there exist some constants $\kappa_2 (\theta ) > \kappa_1 (\theta ) >0$ such that  
\begin{equation}
\kappa_1 (\theta ) \| \phi \| _{{\bf h} (\theta )} \leq \| \mathcal{F}_0 (\theta )^* \phi \| _{\mathcal{B}^* ({\bf Z})} \leq \kappa_2 (\theta ) \| \phi \| _{{\bf h} (\theta )} .
\label{S3_eq_1to1F0*}
\end{equation}
Then $ \mathcal{F}_0 (\theta )^*  $ is one to one.
In particular, the range of $\mathcal{F}_0 (\theta )^* $ is closed.
 For the proof of Theorem \ref{S3_thm_rangeF0*}, we use the following Banach's closed range theorem (see p. 205 of \cite{Yo}).

\begin{theorem}
Let $X_1 ,X_2$ be Banach spaces, and $T$ is a bounded operator in ${\bf B} (X_1 ; X_2 )$.
We denote by $ \langle \cdot , \cdot \rangle $ the pairing between $X_1 $ or $X_2 $ and its dual spaces $(X_1)^*$ or $(X_2)^*$, respectively. 
Then the following assertions are equivalent.
\\
(1) $ \mathrm{Ran} (T) $ is closed. \\
(2) $ \mathrm{Ran} (T^* )$ is closed. \\
(3) $ \mathrm{Ran} (T)= \mathrm{Ker} (T^*)^{\perp} := \{ y\in X_2 \ ; \ \langle y,y^* \rangle =0 , \ y^* \in \mathrm{Ker} (T^*) \} $. \\
(4) $ \mathrm{Ran} (T^* )= \mathrm{Ker} (T)^{\perp} := \{ x^* \in X_1 ^* \ ; \ \langle x,x^* \rangle =0 , \ x \in \mathrm{Ker} (T) \} $.
\label{S3_thm_Banach}
\end{theorem}

Taking $ X_1 = \mathcal{B} ({\bf Z})$, $X_2 = {\bf h} (\theta )$ and $T= \mathcal{F}_0 (\theta )$, we use the assertions (3) and (4).
The equality $ \mathcal{F}_0 (\theta ) \mathcal{B} ({\bf Z}) = {\bf h} (\theta )$ is a direct consequence of the assertion (3), since we have $\mathrm{Ker} \mathcal{F}_0 (\theta )^* =\{ 0\}$ as has been seen in (\ref{S3_eq_1to1F0*}).
For the remaining statement, we have only to prove $(u,f)=0$ when $u\in \mathcal{B}^* ({\bf Z})$, $f\in \mathcal{B} ({\bf Z})$, $(U_0 -e^{i\theta} )u=0$ and $ \mathcal{F}_0  (\theta )f=0$.
Passing to the Fourier series, $(U_0 -e^{i\theta} )u=0 $ is equivalent to $ p(\xi , \theta) \widehat{u} (\xi ) =0$ on ${\bf T}$.
Moreover, this implies $\mathrm{supp} \, \widehat{u} \subset M(\theta )$.
Thus we have as in the case (\ref{S3_eq_L2toH})
\begin{gather*}
\begin{split}
&(u,f)\\
 & = \sum_{j=1}^2 \int _{J_{\gamma ,j}} \sum _{\xi (\omega )\in M(\omega )} (P_j (\xi (\omega )) \widehat{u} (\xi (\omega )) , P_j (\xi ( \omega )) \widehat{f} (\xi (\omega )) )_{{\bf C}^2}  \frac{d\xi}{d\omega} (\omega )  d\omega \\
&= \sum _{j=1}^2 \sum _{\xi (\theta )\in M(\theta )} (P_j (\xi (\theta )) \widehat{u} (\xi (\theta )) , P_j (\xi ( \theta )) \widehat{f} (\xi (\theta )) )_{{\bf C}^2}  \frac{d\xi}{d\omega} (\theta )  .
\end{split}
\end{gather*}
Then $ \mathcal{F}_0 (\theta )f=0$ implies $ (u,f)=0$.
This proves the lemma.
\qed


\section{Position-dependent QW}

\subsection{Generalized eigenfunction with radiation conditions}

Let us turn to the operator $U$ for the position-dependent quantum walk.
For the beginning, we consider generalized eigenfunction $u\in \mathcal{B}^* ({\bf Z})$ satisfying the equation
\begin{equation}
(U-e^{i\theta} ) u=0 , \quad \theta \in J_{\gamma} \setminus J_{\gamma, \mathcal{T}} .
\label{S4_eq_evprob}
\end{equation}
In the following, we put
\begin{equation}
V=U-U_0 .
\label{S4_def_V}
\end{equation}

We will construct generalized eigenfunctions satisfying the equation (\ref{S4_eq_evprob}) in $\mathcal{B}^* ({\bf Z} )$ later.
Here we prove the absence of non-trivial solutions satisfying the condition (\ref{S3_eq_radconditionlattice}).

\begin{lemma}
Let $\theta \in J_{\gamma} \setminus J_{\gamma , \mathcal{T} } $.
Suppose $u_{\pm} \in \mathcal{B}^* ({\bf Z})$ satisfy the equation (\ref{S4_eq_evprob}) and the condition 
$$
WF^* ( P_1 \mathcal{U} u_{\pm} )= \{ (\xi (\theta ), \mp 1)\} , \quad WF^* ( P_2 \mathcal{U} u_{\pm} )= \{ (\xi (\theta ), \pm 1)\} ,
$$
respectively.
Then we have $ u_{\pm} \in \mathcal{B}_0^* ({\bf Z})$.
\label{S4_lem_efB0}
\end{lemma}

Proof.
In view of the assumption (A-1), we note that $ |(Vu)(x)| \leq M e^{-\epsilon_0 \langle x\rangle} $ for constants $ \epsilon_0 ,M>0$.

In the following, we discuss $ u_+ $.
Letting $ f_+ = -Vu_+ $, we have $ (U_0 -e^{i\theta} ) u_+ = f_+ $.
We also have $ (U_0^* - e^{-i\theta} ) u_+= f_+^* $ where $f_+^* = -V^* u_+ $, $V^* = U^* - U_0^*$. 
Note that $f_+  = -e^{i\theta } U_0 f_+^* $.
The assertion (4) in Theorem \ref{S3_thm_LAPlattice} and (\ref{S3_eq_WFadjoint}) imply $ u_+ = R_0 (\theta +i0 )f_+ = R_0 (\theta -i0 )^* f_+^* $.

Now we compute $ ( R_0 (\theta +i0) f_+ - R_0 (\theta -i0 ) f_+ , f_+ )$.
In fact, we have
\begin{gather*}
\begin{split}
&( R_0 (\theta +i0) f_+ - R_0 (\theta -i0 ) f_+ , f_+ ) \\
=& \, - (u_+ , V u_+ ) -   ( U_0 f_+^* , R_0 (\theta -i0 )^* U_0 f_+^* ) \\
= & \, -( u_+ , V u_+ ) - (U_0 f_+^* , U_0 u_+ )  \\
= & \, -( u_+ , V u_+ ) + ( V^* u_+ , u_+ )=0 .
\end{split}
\end{gather*}
On the other hand, it follows from the assertion (5) in Theorem \ref{S3_thm_LAPlattice} that 
$$
( R_0 (\theta +i0 ) f_+ - R_0 (\theta -i0) f_+ , f_+ )= 2\pi e^{-i\theta} \| \mathcal{F}_0 (\theta ) f_+ \| ^2 _{{\bf h} (\theta )} .
$$
Thus we have $ \mathcal{F}_0 (\theta ) f_+ =0$.
We apply again the assertion (4) in Theorem \ref{S3_thm_LAPlattice} and we obtain $u_+ \in \mathcal{B}_0^* ({\bf Z})$.

For $u_-$, the proof is similar.
We have proven the lemma.
\qed

\begin{theorem}
Let $\theta \in J_{\gamma} \setminus J_{\gamma , \mathcal{T}} $ and $ u _{\pm} \in \mathcal{B} ^*_0 ({\bf Z})$ be as in Lemma \ref{S4_lem_efB0}. 
Then we have $ u_{\pm} =0 $.

\label{S4_thm_absenceofef}
\end{theorem}

Proof.
In view of Theorem \ref{S2_thm_embev}, we have only to show $ u_{\pm} \in \mathcal{H} $.
We shall prove for $ u_+ $.
Passing to the Fourier series, the equation $(U-e^{i\theta} )u_+ =0$ is rewritten by
\begin{equation}
(\widehat{U}_0 (\xi ) - e^{i\theta} ) \widehat{u}_+ (\xi ) = \widehat{f}_+ (\xi ) \quad \text{on} \quad {\bf T} , 
\label{S4_eq_eigeneqtorus}
\end{equation}
where $\widehat{f} _+ =  -\mathcal{U} (Vu_+) $.
Since we have $ |(V u_+)(x) |_{{\bf C}^2} \leq M e^{-\epsilon_0 \langle x\rangle} $, the Payley-Wiener theorem implies that $\widehat{f}_+ $ extends to an analytic function in a neighborhood of ${\bf T}$ (see Theorem 6.1 in \cite{Ves}).
We multiply the equation (\ref{S4_eq_eigeneqtorus}) by the cofactor matrix of $ \widehat{U}_0 (\xi ) -e^{i\theta} $. 
Then the left-hand side of the equation (\ref{S4_eq_eigeneqtorus}) is diagonalized as 
\begin{equation}
p(\xi , \theta ) \widehat{u}_+ (\xi ) = \widehat{g}_+(\xi ) ,
\label{S4_eq_eigeneqtorus2}
\end{equation}
where $\widehat{g} _+ (\xi ) $ extends to an analytic function in a neighborhood of ${\bf T} $.

For $\xi (\theta )\in M(\theta )$, we take $\chi \in C^{\infty} ({\bf T})$ such that $\chi (\xi (\theta ))=1 $ with small support.
Since we have assumed $\theta \in  J_{\gamma} \setminus J_{\gamma , \mathcal{T}} $, we have $ M(\theta )= M_{reg} (\theta )$ so that we can make the change of variable $\xi \mapsto \eta $ as in \S 3 in a small neighborhood of $\xi (\theta )$.
Letting $ \widehat{u}_{\chi} = \chi \widehat{u}_+ $ and $ \widehat{g}_{\chi } = \chi \widehat{g}_+ $, we rewrite the equation (\ref{S4_eq_eigeneqtorus2}) as 
$$
\eta \widehat{u} _{\chi} = -\frac{1}{2p} e^{ -i (\theta + \gamma /2)} \widehat{g} _{\chi } . 
$$
Moreover, by using the Fourier transformation, we have 
\begin{equation}
\frac{d \widetilde{u}_{\chi}}{d s} = \frac{i}{2p} e^{-i (\theta + \gamma /2)} \widetilde{g}_{\chi} .
\label{S4_eq_eigeneqtorus3}
\end{equation}
Integrating this equation, we have
$$
\widetilde{u} _{\chi} (s)= \frac{i}{2p} e^{-i (\theta + \gamma /2)}  \int_0^s \widetilde{g} _{\chi} (\tau) d\tau + \widetilde{u}_{\chi} (0) .
$$
Since $\widetilde{g}_{\chi} $ is rapidly decreasing, the limit 
\begin{equation}
\lim_{s \to  \pm \infty} \widetilde{u} _{\chi} (s)= \frac{i}{2p} e^{-i (\theta + \gamma /2)}  \int_0^{\pm \infty} \widetilde{g} _{\chi} (\tau) d\tau + \widetilde{u}_{\chi} (0)  ,
\label{S4_eq_limit}
\end{equation}
exists.

It follows from Lemma \ref{S4_lem_efB0} that $ \widetilde{u} _{\chi } \in \mathcal{B}_0^* ({\bf R})$.
Thus we have
$$
\frac{1}{R} \int _{R/2<|s|<R} | \widetilde{u}_{ \chi } (s) |^2 _{{\bf C}^2} ds \leq \frac{1}{R} \int _{|s|<R} | \widetilde{u}_{ \chi } (s) |^2 _{{\bf C}^2} ds  \to 0 ,
$$
as $R\to \infty $.
This implies $ \liminf _{|s|\to \infty} |\widetilde{u}_{\chi} (s)| _{{\bf C}^2}=0 $.
We have therefore by (\ref{S4_eq_limit})
$$
\widetilde{u}_{\chi} (0)= -\frac{i}{2p} e^{-i (\theta + \gamma /2)}  \int _0^{\infty} \widetilde{g} _{\chi} (\tau ) d\tau = \frac{i}{2p} e^{-i (\theta + \gamma /2)} \int _{-\infty}^0 \widetilde{g}_{\chi} (\tau ) d\tau .
$$
The solution $\widetilde{u} _{\chi } $ is represented by 
$$
\widetilde{u} _{\chi } (s)= \frac{i}{2p} e^{-i (\theta + \gamma /2)} \int_s^{\infty} \widetilde{g}_{\chi} (\tau ) d\tau , \quad s\geq 0 , 
$$
and
$$
\widetilde{u} _{\chi } (s)= \frac{i}{2p} e^{-i (\theta + \gamma /2)} \int^s_{-\infty} \widetilde{g}_{\chi} (\tau ) d\tau , \quad s\leq 0 .
$$
Then $ \widetilde{u} _{\chi } (s)$ is also rapidly decreasing as $|s| \to \infty $.
This implies $\widehat{u}_{\chi} \in C^{\infty } ({\bf T})$.
Obviously, $\widehat{u}_+$ is smooth outside any small neighborhood of $\xi (\theta )$.
Thus we obtain $\widehat{u} _+ \in C^{\infty} ({\bf T})$. 
We have proven $u_+ \in \mathcal{H} $.
\qed


\subsection{Resolvent estimate and spectral representation}

We put 
\begin{equation}
R(z)= (U-e^{iz} )^{-1} , \quad z\in {\bf C} \setminus {\bf R} .
\label{S4_eq_resolvent}
\end{equation}
The well-known resolvent equation holds :
\begin{equation}
R(z)= R_0 (z) - R_0 (z) V R(z) .
\label{S4_eq_resolventeq_pert}
\end{equation}
Now we derive the limiting absorption principle for $R(z)$.

\begin{theorem}
Let $J \subset J_{\gamma} \setminus J_{\gamma , \mathcal{T}} $ be a compact interval and $\theta$ vary in $J$.
\\
(1) There is a constant $ \kappa >0$ such that 
$$
\| R(z)f \| _{\mathcal{B}^* ({\bf Z})} \leq \kappa \| f\| _{\mathcal{B} ({\bf Z})} , \quad \mathrm{Re} \, z \in J , \quad \mathrm{Im} \, z \not= 0.
$$
(2) There exists a limit $ R(\theta \pm i0 ) := \lim _{\epsilon \downarrow 0 } R(\theta -i \log (1\mp \epsilon )) $ in the weak $*$ sense.
Moreover, we have $R(\theta \pm i0) \in {\bf B} (\mathcal{B} ({\bf Z}) ; \mathcal{B}^* ({\bf Z} ))$ with 
$$
\| R( \theta \pm i0 )f\| _{\mathcal{B}^* ({\bf Z})} \leq \kappa \| f\| _{\mathcal{B} ({\bf Z} )} ,
$$
for a constant $ \kappa >0 $. \\
(3) The mapping 
$$
J_{\gamma} \setminus J_{\gamma , \mathcal{T}} \ni \theta \mapsto ( R(\theta \pm i0 ) f,g) , \quad f,g\in \mathcal{B} ({\bf Z}) ,
$$
is continuous.

\label{S4_thm_LAP}
\end{theorem}

Proof.
In view of (\ref{S4_eq_resolventeq_pert}), we note the inequality
\begin{equation}
\| R(z) \| _{\mathcal{B}^* ({\bf Z})} \leq \kappa \left( \| f\| _{\mathcal{B} ({\bf Z})} + \| V R(z)f\| _{\mathcal{B} ({\bf Z})} \right) .
\label{S4_eq_resolventestimate11}
\end{equation}

First let us show the assertion (1).
Suppose that the assertion (1) does not hold.
Then there exist sequences $\{ f_j \} _{j=1,2,\ldots } \subset \mathcal{B} ({\bf Z})$ and $\{ z_j \} _{j=1,2,\ldots} \subset {\bf C} \setminus {\bf R} $ such that 
$$
\| R(z_j ) f_j \| _{\mathcal{B}^* ({\bf Z} )} =1 , \quad \| f_j \| _{\mathcal{B} ({\bf Z}) } \to 0 , 
$$
as $j\to \infty $.
In the following, we assume $z_j \to \theta +i0$ without loss of generality.
We put $u_j = R(z_j ) f_j $. 
For $ \sigma > 1/2 $, the embedding $\mathcal{B}^* ({\bf Z} ) \subset \ell^{2,-\sigma} ({\bf Z})$ is compact.
Then we can take a subsequence of $\{ u_j \} $ such that it converges to $u \in \ell^{2,-\sigma} ({\bf Z})$.
We denote by $\{ u_j \}$ this subsequence again.
The inequality (\ref{S4_eq_resolventestimate11}) implies $u\in \mathcal{B}^* ({\bf Z} )$ with $ \| u\| _{\mathcal{B}^* ({\bf Z})} \leq \kappa \| V u\| _{\mathcal{B} ({\bf Z})} $ so that $\| u\| _{\mathcal{B}^* ({\bf Z})} =1 $ by the assumption.
On the other hand, $u$ satisfies $ (U -e^{i\theta} )u= 0$ and 
$$
u= - R_0 ( \theta +i0 )Vu .
$$
Thus we can apply Lemma \ref{S4_lem_efB0} and Theorem \ref{S4_thm_absenceofef} to $u$ i.e. we have $u=0$.
This is a contradiction.

Let us turn to the assertion (2).
We take a sequence $z_j = \theta -i \log (1-\epsilon_j )$ with $\epsilon_j \downarrow 0$.
For $f\in \mathcal{B} ({\bf Z})$, we put $ u_j = R(\theta -i\log (1-\epsilon_j )) f$.
In view of the assertion (1), there exists a subsequence of $\{ u_j \}$, which is denoted by $\{ u_j \}$ again, such that $u_j \to u$ in $\ell^{2,-\sigma} ({\bf Z} )$.
It follows $u\in \mathcal{B}^* ({\bf Z} )$ from the resolvent equation 
$$
u_j = R_0 (z_j )f - R_0 (z_j ) Vu_j \to R_0 (\theta +i0 )f - R_0 (\theta +i0 ) Vu ,
$$
as $j\to \infty $.

We prove that $\{ u_j \} $ itself converges to $u$ in $\ell^{2,-\sigma} ({\bf Z} )$.
Assume that there exists a subsequence of $\{ u_j \}$ such that it does not converges to $u$.
We denote by $\{ u_{j_k} \} $ this subsequence. 
Then we can take a constant $\delta >0$ such that 
$$
\| u_{j_k} -u \| _{\ell^{2,-\sigma} ({\bf Z})} \geq \delta ,
$$
for any $k =1,2, \ldots $.
The subsequence $\{ u_{j_k} \}$ has a sub-subsequence such that it converges to $u' $ in $\ell^{2,-\sigma} ({\bf Z} )$.
Moreover, we can prove $u'\in \mathcal{B}^* ({\bf Z})$ and that $u'$ satisfies 
$$
u'= R_0 (\theta +i0 )f - R_0 (\theta +i0 ) Vu ' .
$$
Then $u-u' $ satisfies $ (U-e^{i\theta} )(u-u')=0$ and 
$$
u-u'= -R_0 (\theta +i0 )V(u-u') .
$$
Lemma \ref{S4_lem_efB0} and Theorem \ref{S4_thm_absenceofef} imply $u = u'$ which is a contradiction.

Let us show that the limit $u_j \to u:= R(\theta +i0)f$ exists in the weak $*$ sense. 
For $f,g \in \ell^{2,\sigma} ({\bf Z} )$, the limit $ u_j \to u$ in $\ell^{2,-\sigma} ({\bf Z})$ implies $( R(\theta -i\log (1-\epsilon )) f,g) \to (R(\theta +i0)f,g)$. 
By using the assertion (1), we have
\begin{equation}
|( R(\theta +i0)f,g)| = \lim _{\epsilon \downarrow 0 } | (R(\theta -i\log (1-\epsilon )) f,g)| \leq \kappa \| f\| _{\mathcal{B} ({\bf Z}) } \| g\| _{\mathcal{B} ({\bf Z} )} .
\label{S4_eq_weak*11}
\end{equation}
The space $\ell^{2,\sigma} ({\bf Z})$ is dense in $\mathcal{B} ({\bf Z} )$.
Thus the above weak $*$ limit exists for any $f,g\in \mathcal{B} ({\bf Z} )$.

The assertion (3) follows from the inequality (\ref{S4_eq_weak*11}) and the denseness of $\ell^{2,\sigma} ({\bf Z})$ in $\mathcal{B} ({\bf Z} )$.

For $R(\theta -i0)$, the assertions (2) and (3) can be proven by the similar way.
\qed

\medskip

We define
\begin{equation}
\mathcal{F} _{\pm} (\theta )= \mathcal{F}_0 (\theta ) (1-VR(\theta \pm i0 ))  , \quad \theta \in J_{\gamma} \setminus J_{\gamma , \mathcal{T}} .
\label{S4_def_fpm}
\end{equation}

\begin{lemma}
For $ \theta \in J_{\gamma} \setminus J_{\gamma , \mathcal{T}} $, we have
$$
(R(\theta +i0 )f -R(\theta -i0 )f ,g)= 2\pi e^{-i\theta} ( \mathcal{F} _{\pm} (\theta )f ,\mathcal{F}_{\pm} (\theta ) g )_{{\bf h} (\theta )} ,
$$
for $f,g \in \mathcal{B} ({\bf Z} )$, and $ \mathcal{F}_{\pm} (\theta ) \in {\bf B} (\mathcal{B} ({\bf Z}) ; {\bf h} (\theta ))$ with the estimate
$$
\| \mathcal{F}_{\pm} (\theta )f \|_{{\bf h} (\theta )} \leq \kappa \| f\| _{\mathcal{B} ({\bf Z})} , \quad f\in \mathcal{B} ({\bf Z} ),
$$
for a constant $\kappa >0$.
\label{S4_lem_stone}
\end{lemma}

Proof.
By the resolvent equation, we have 
\begin{gather*}
\begin{split}
&R ( \theta - i \log ( 1-\epsilon )) - R( \theta -i\log (1+\epsilon )) \\
&= -2\epsilon e^{i\theta} R( \theta - i\log (1-\epsilon )) R( \theta -i \log (1+\epsilon )) \\
&= -2\epsilon e^{i\theta} R( \theta - i\log (1+\epsilon )) R( \theta -i \log (1-\epsilon )) ,
\end{split}
\end{gather*}
for $\epsilon > 0$.
Thus, for $f,g\in \mathcal{B} ({\bf Z})$, we have 
\begin{gather*}
\begin{split}
&( R(\theta +i0) f - R(\theta -i0) f,g) \\
&= ( (R_0 (\theta +i0)- R_0 (\theta -i0 ))(1-VR (\theta \pm i0) )f, (1-V^* R(\theta \mp i0)^* ) g) \\
&= 2\pi e^{-i \theta } ( \mathcal{F}_0 (\theta ) (1-V R(\theta \pm i0 )) f, \mathcal{F}_0 (\theta ) (1-V^* R(\theta \mp i0)^*) g ) _{{\bf h} (\theta )} .
\end{split}
\end{gather*}
The equality $ R(\theta \pm i0)^* = -e ^{i\theta} U R(\theta \mp i0)$ implies 
\begin{gather*}
\begin{split}
\mathcal{F}_0 (\theta ) (1-V^* R( \theta \mp i0 )^* ) 
&= \mathcal{F}_0 (\theta )(1-e^{i \theta} U_0^* (U-U_0 ) R(\theta \pm i0 )) \\
&= \mathcal{F}_0 (\theta )(1-V R(\theta \pm i0 )) . 
\end{split}
\end{gather*}
Therefore, we obtain the lemma.
\qed

\medskip

Let us define $ (\mathcal{F}_{\pm} f )(\theta )=  \mathcal{F}_{\pm} (\theta )f$ for $f\in \mathcal{B} ({\bf Z} )$. 
In order to derive an extension of $\mathcal{F}_{\pm}$ to $\mathcal{H}_{ac} (U)$, we show an analogue of Stone's formula.

\begin{lemma}
We have 
\begin{gather*}
\begin{split}
&\lim _{\epsilon \downarrow 0} \int_{\theta_1}^{\theta_2} \frac{e^{i\theta}}{2\pi} ( (R( \theta -i\log (1-\epsilon )) f -R(\theta -i\log (1+\epsilon )) f ,g) d\theta \\
&= E_U ((\theta_1 , \theta _2 )) ,
\end{split}
\end{gather*}
for $ \theta_1 < \theta_2  $ and $ \theta _1 , \theta_2 \in J_{\gamma} \setminus J_{\gamma , \mathcal{T}} $.
\label{S4_lem_stone2}
\end{lemma}

Proof.
In view of the spectral decomposition of $U$, we have 
\begin{gather*}
\begin{split}
& \int_{\theta_1}^{\theta_2} \frac{e^{i\theta}}{2\pi} ( (R( \theta -i\log (1-\epsilon )) f -R(\theta -i\log (1+\epsilon )) f ,g) d\theta \\
&= \frac{1}{2\pi} \int _{\sigma (U)}  (S_{+,\epsilon } (\omega )- S_{-,\epsilon} (\omega )) d(E_U (\omega )f,g) _{\mathcal{H}} ,
\end{split}
\end{gather*}
where $ S_{\pm ,\epsilon } (\omega )$ is defined by (\ref{A_eq_stone}).
Note that $ \sigma_p (U) \cap (\sigma_{ess} (U) \setminus \mathcal{T})=\emptyset $.
Then we have $ E_U (\{ \theta_1 \})= E_U (\{ \theta_2 \} )=0$.
Thus Lemma \ref{A_lem_stoneformula} implies 
\begin{gather*}
\begin{split}
& \lim_{\epsilon \downarrow 0} \int_{\theta_1}^{\theta_2} \frac{e^{i\theta}}{2\pi} ( (R( \theta -i\log (1-\epsilon )) f -R(\theta -i\log (1+\epsilon )) f ,g) d\theta \\
& = \int_{\theta_1}^{\theta_2} d(E_U (\omega )f,g)_{\mathcal{H}}.
\end{split}
\end{gather*}
This equality shows the lemma.
\qed

\medskip

It is well-known for self-adjoint operators that the limiting absorption principle and Stone's formula imply the absence of the singular continuous spectrum.
We can also show the absence of singular continuous spectrum of $U$ as follows.

\begin{lemma}
We have $\sigma_{sc} (U)=\emptyset $.
\label{S4_lem_scspec}
\end{lemma}

Proof.
Let $\theta \in J_{\gamma} \setminus J_{\gamma , \mathcal{T}} $ and $ I=( \theta - \epsilon , \theta + \epsilon )$ for small $\epsilon >0$.
It follows from Lemma \ref{S4_lem_stone2} that 
$$
(E_U ((\theta -\epsilon , \theta '  )) u,u) _{\mathcal{H}} = \int _{\theta -\epsilon}^{\theta '} \frac{e^{i\omega} }{2\pi } (R( \omega +i0)u- R(\omega -i0)u,u) d\omega ,
$$
for $\theta '\in I$ and $u\in \mathcal{B} ({\bf Z})$.
Then the continuity of $(R( \omega +i0)u- R(\omega -i0)u,u)$ implies $ E_U (I) u \in \mathcal{H}_{ac} (U)$.
Since $ \mathcal{B} ({\bf Z})$ is dense in $ \mathcal{H}$ and $\mathcal{H}_{ac} (U)$ is closed, we have $E_U (I) \mathcal{H} \subset \mathcal{H}_{ac} (U)$.

Suppose $\theta ' \in \sigma_{sc} (U) \cap I$. 
Then there exists a sequence $\{ u_j \} \subset \mathcal{H}_{sc} (U)$ such that $ \| u_j \| _{\mathcal{H}} =1$ and $\| (U-e^{i\theta '} )u_j \| _{\mathcal{H}} \to 0$ as $j\to \infty $.
However $E_U (I) u_j =0$ implies 
\begin{gather*}
\begin{split}
\| (U-e^{i\theta '} )u_j \|^2 _{\mathcal{H}} = & \int _{|\omega -\theta' | >\epsilon } | e^{i\omega} - e^{i\theta '} |^2 d (E_U (\omega )u_j , u_j ) _{\mathcal{H}}  \\
\geq & \epsilon^2 \| u_j \|^2_{\mathcal{H}} = \epsilon^2 >0 ,
\end{split}
\end{gather*}
for any $j$.
This is a contradiction.
\qed

\medskip

Now we can state the spectral representation of $U$ as a partial isometry between $ \mathcal{H}_{ac} (U) $ and $ {\bf H} $ as follows.

\begin{theorem}
(1) The operator $\mathcal{F}_{\pm} $ can be extended to a partial isometry with initial set $ \mathcal{H}_{ac} (U) $ and final set ${\bf H}$. \\
(2) We have $( \mathcal{F}_{\pm} U f)(\theta )=e^{i\theta} (\mathcal{F}_{\pm} f)(\theta )$ and $( \mathcal{F}_{\pm} U ^* f)(\theta )=e^{-i\theta} (\mathcal{F}_{\pm} f)(\theta )$ for $f\in \mathcal{H}_{ac} (U)$ and $\theta \in J_{\gamma} \setminus J_{\gamma , \mathcal{T}}$. \\
(3) For $\theta \in J_{\gamma} \setminus J_{\gamma , \mathcal{T}} $ and $\phi \in {\bf h} (\theta )$, $\mathcal{F}_{\pm} (\theta )^* \phi \in \mathcal{B}^* ({\bf Z})$ satisfies $(U-e^{i\theta}) \mathcal{F}_{\pm } (\theta )^* \phi =0$.
\label{S4_thm_extensionFpm}
\end{theorem}

Proof.
It follows from Lemmas \ref{S4_lem_stone} and \ref{S4_lem_stone2} that 
$$
\int_{\theta_1}^{\theta_2} ( \mathcal{F}_{\pm} (\theta) f,\mathcal{F}_{\pm} (\theta_2 )g)_{{\bf h} (\theta)} d\theta = (E_U ((\theta_1 , \theta_2 )) f,g) _{\mathcal{H} } ,
$$
for $\theta_1 < \theta_2 $, $\theta_1 , \theta_2 \in J_{\gamma} \setminus J_{\gamma , \mathcal{T}}$.
Then we can see the assertion (1).

It follows $\mathcal{F}_{\pm} (\theta ) U f = e^{i\theta} \mathcal{F}_{\pm} (\theta )f$ for $f\in \mathcal{B} ({\bf Z}) \cap \mathcal{H}_{ac} (U)$ from a direct computation.
Let us show $ \mathcal{F}_{\pm} (\theta ) U^* f = e^{-i\theta} \mathcal{F}_{\pm} (\theta )f$ for $ f\in \mathcal{B} ({\bf Z} ) \cap \mathcal{H}_{ac} (U) $. 
Note that $ R( \theta \pm i0 ) U^* f= e^{-i\theta} R(\theta \pm i0) f - e^{-i\theta} U^* f$ for $ f\in \mathcal{B} ({\bf Z} ) \cap \mathcal{H}_{ac} (U) $. 
Thus we have
\begin{gather*}
\begin{split}
\mathcal{F}_{\pm} (\theta )U^* f &=   e^{-i\theta} \mathcal{F}_{0} (\theta )f + \mathcal{F}_0 (\theta ) V^* f - \mathcal{F}_0 (\theta )VR (\theta \pm i0)U^* f \\
&= e^{-i\theta} \mathcal{F}_{\pm} (\theta )f + \mathcal{F}_0 (\theta ) \left( V^* f + e^{-i\theta}   VU^* f \right) .
\end{split}
\end{gather*}
Moreover, we can see
$$
\mathcal{F}_0 (\theta ) \left( V^* f + e^{-i\theta}   VU^* f \right) = 0.
$$
Then we have $ \mathcal{F}_{\pm} (\theta )U^* f= e^{-i\theta} \mathcal{F}_{\pm} (\theta )f$.
This is a proof of the assertion (2).

Taking the adjoint, we also have the assertion (3).
\qed

\medskip


\section{Wave operator and scattering matrix}
\subsection{Scattering matrix}

We consider the wave operators $W_{\pm} $ defined by 
\begin{equation}
W_{\pm} = {\mathop{{\rm s-lim}}_{t\to\pm \infty}} \, U^{-t} U_0^t , \quad t\in {\bf Z} .
\label{S5_def_wave operators}
\end{equation}
The existence and the completeness of $ W_{\pm} $ have been proven by \cite{Su} and \cite{RST} (under an assumption which is weaker than that of the present paper) as follows.
\begin{theorem}
The wave operators exist and are complete i.e. $ \mathrm{Ran} W_{\pm} = \mathcal{H}_{ac} (U)$.
As a consequence, for any $ \phi \in \mathcal{H} _{ac} (U) $, there exist $ \psi _{\pm} \in \mathcal{H}$ such that $ \| U_0^t \phi - U^t \psi _{\pm} \| _{\mathcal{H}} \to 0$ as $t\to \pm \infty $.
The wave operators are unitary on $ \mathcal{H}$ and we have $ (W_{\pm} )^* = (W_{\pm} )^{-1} $.
\label{S5_thm_WO}
\end{theorem}

Thus the following sums converge in ${\bf B} (\mathcal{H} )$ :
\begin{gather}
W_+ = 1 + \sum _{t=0}^{\infty} U^{-t} V^* U_0 ^{t+1}  , \label{S5_eq_W+} \\
W_- = 1 - \sum _{t=-1}^{-\infty} U^{-t} V^* U_0 ^{t+1}  . \label{S5_eq_W-}
\end{gather}

In the following, we relate the wave operators $ W_{\pm} $ with the resolvent operators $ R (\theta  \pm i0 )$ and $ R_0 (\theta \pm i0)$.
Finally, we will construct a explicit formula of the scattering operator :
\begin{equation}
{\bf S} := (W_+)^* W_- .
\label{S5_def_Soperator}
\end{equation}

\begin{lemma}
We have
\begin{gather}
\sum _{t=0}^{\infty} e^{-i t \theta } U^{t} =- e^{i\theta} R(\theta -i0) , \label{S5_eq_WO1} \\
\sum_{t=-\infty}^{\infty} e^{-it\theta} U^t_0 = e^{i\theta} (R_0 (\theta +i0) -R_0 (\theta -i0)) =2\pi \mathcal{F}_0 (\theta )^* \mathcal{F}_0 (\theta ) , \label{S5_eq_WO2} 
\end{gather}
in the weak $*$ sense for $ \theta \in J_{\gamma} \setminus J_{\gamma , \mathcal{T}} $.

\label{S5_lem_WO1}
\end{lemma}

Proof.
Let $ F_{\epsilon} (t) = e^{-\epsilon t} e^{-it \theta} U^{t} f$ for $\epsilon >0$ and $f\in \mathcal{H} $.
Obviously, we have 
$$
\sum _{t=0}^{\infty} \| F_{\epsilon} (t)\| ^n _{\mathcal{H}} \leq \sum _{t=0} ^{\infty} e^{-\epsilon tn} \| f\| ^n _{\mathcal{H}} <\infty ,
$$
for $ n=1,2 $.
In view of the equality
$$
(U^{t} f,g) _{\mathcal{H}} = \int _0^{2\pi} e^{-it \omega } d(E_U (\omega ) f,g) _{\mathcal{H}} , \quad f,g \in \mathcal{B} ({\bf Z} ) ,
$$
we obtain
\begin{gather*}
\begin{split}
\left( \sum_{t=0}^{\infty} F_{\epsilon} (t) ,g \right) _{\mathcal{H}} &= \sum_{t=0}^{\infty} e^{-\epsilon t} e^{-it \theta } (U^{t} f,g) _{\mathcal{H}} \\
&= \int_0^{2\pi} \sum_{t=0}^{\infty} e^{-\epsilon t} e^{-it \theta} e^{it \omega} d(E_U (\omega )f,g) _{\mathcal{H}} \\
&= -e^{\epsilon} e^{i\theta} \int_0^{2\pi} \frac{1}{e^{i\omega} - e^{i\theta} e^{\epsilon} } d(E_U (\omega ) f,g )_{\mathcal{H}} .
\end{split}
\end{gather*}
Tending $\epsilon \downarrow 0$, we obtain (\ref{S5_eq_WO1}).

Let $ F'_{\epsilon} (t)= e^{-\epsilon |t|} e^{-it \theta} U_0^t f$ for $ \epsilon >0$ and $ f\in \mathcal{H} $.
As above, we have $ \sum _{t=-\infty}^{\infty} \| F'_{\epsilon} (t) \| _{\mathcal{H}}^n <\infty $ for $n=1,2$.
We also have
$$
\left( \sum _{t=-\infty}^{\infty} F'_{\epsilon} (t) ,g \right)_{\mathcal{H}} = \int_0^{2\pi} \sum _{t=-\infty}^{\infty} e^{-\epsilon |t|} e^{-it\theta} e^{it\omega} d(E_{U_0} (\omega ) f,g) _{\mathcal{H}} .
$$
Note that 
$$
\sum _{t=-\infty}^{\infty} e^{-\epsilon |t|} e^{-it\theta} e^{it\omega} = \frac{e^{-\epsilon} e^{i\theta}}{e^{i\omega} - e^{i\theta} e^{-\epsilon} } - \frac{e^{\epsilon} e^{i\theta} }{e^{i\omega} - e^{i\theta} e^{\epsilon} } .
$$
Hence, tending $ \epsilon \downarrow 0 $, we have the first equality of (\ref{S5_eq_WO2}).
The second equality of (\ref{S5_eq_WO2}) is a direct consequence of the assertion (5) of Theorem \ref{S3_thm_LAPlattice}. 
\qed

\medskip

Now we define the \textit{S-matrix} $\widehat{{\bf S}} (\theta )$ by
\begin{equation}
\widehat{{\bf S}} = \mathcal{F}_0 {\bf S} \mathcal{F}_0^* = \int_{J_{\gamma}} \oplus \widehat{{\bf S}} (\theta ) d\theta .
\label{S5_eq_Smatrix}
\end{equation}
Then $\widehat{{\bf S}} (\theta )$ is an operator on ${\bf h} (\theta )$.

\begin{theorem}
(1) For $f\in {\bf H}$, we have $ (\widehat{{\bf S}} f) (\theta ) = \widehat{{\bf S}} (\theta ) f(\theta )$ for every $ \theta \in J_{\gamma} \setminus J_{\gamma , \mathcal{T}} $. \\
(2) $\widehat{{\bf S}} (\theta )$ is unitary on ${\bf h} (\theta )$ for $\theta \in J_{\gamma} \setminus J_{\gamma , \mathcal{T}} $. \\
(3) We have 
\begin{equation}
\widehat{{\bf S}} (\theta )=1 -2\pi e^{i\theta} A(\theta ),
\label{S5_eq_Smatrix_ampli}
\end{equation}
where
$$
A (\theta )= \mathcal{F}_0 (\theta ) \left( V^* - VR(\theta -i0 )V^*  \right) \mathcal{F}_0 (\theta )^* . 
$$
for $ \theta \in J_{\gamma} \setminus J_{\gamma , \mathcal{T}} $.
\label{S5_thm_Smatrix}
\end{theorem}

Proof.
Noting the equalities (\ref{S5_eq_W+}) and (\ref{S5_eq_W-}), let us compute $ (({\bf S}-1) E_{U_0} (I)f, E_{U_0} (I) g) _{\mathcal{H}} $ for any compact interval $ I \subset J_{\gamma} \setminus J_{\gamma , \mathcal{T}} $ and $ f,g\in \mathcal{H} $.
To begin with, we take $ f,g \in \mathcal{B} ({\bf Z} )$.
By the definition, we have $ {\bf S} -1 = (W_+)^* (W_- - W_+ )$ so that 
\begin{gather*}
\begin{split}
&(({\bf S}-1) E_{U_0} (I)f, E_{U_0} (I) g) _{\mathcal{H}} \\
&=( (W_- - W_+ ) E_{U_0} (I) f , W_+ E_{U_0} (I) g )_{\mathcal{H}} \\
&= - \sum_{t=-\infty}^{\infty} (V^* U_0^{t+1} E_{U_0} (I) f, U_0^t E_{U_0} (I) g)_{\mathcal{H}} \\
&\quad \ \ \ \ - \sum_{t=-\infty}^{\infty} \sum_{\tau =0}^{\infty} (V^* U_0^{t+1} E_{U_0} (I) f, U^{-\tau} V^* U_0^{t+\tau +1} E_{U_0} (I)g )_{\mathcal{H}} \\
&=: - {\bf S}_1 - {\bf S}_2  .
\end{split}
\end{gather*}
Here we have used the equality $ U^{\tau} W_{\pm} = W_{\pm} U_0^{\tau} $. 

First we shall compute ${\bf S}_1 $ as follows :
$$
{\bf S}_1 = \sum _{t=-\infty}^{\infty} \int_I ( \mathcal{F}_0 (\theta ) V^* U_0^{t+1} E_{U_0} (I) f, e^{it\theta} \mathcal{F}_0 (\theta ) g) _{{\bf h} (\theta )} d\theta . 
$$
In view of the formula (\ref{S5_eq_WO2}), we have $ \sum_{t=-\infty}^{\infty} e^{-it \theta} U_0^{t+1} =2\pi e^{i\theta} \mathcal{F}_0 (\theta )^* \mathcal{F}_0 (\theta )$, inserting $e^{-\epsilon |t+1|} $ and letting $\epsilon \downarrow 0$.
Then ${\bf S}_1 $ can be rewritten as 
\begin{equation}
{\bf S}_1 = 2\pi \int_I e^{i\theta} ( \mathcal{F}_0 (\theta ) V^* \mathcal{F}_0 (\theta )^* \mathcal{F}_0 (\theta )f, \mathcal{F}_0 (\theta )g)_{{\bf h} (\theta )} d\theta .
\label{S5_eq_Smatrix11}
\end{equation}

Let us turn to ${\bf S}_2 $. 
Obviously, we have
\begin{gather*}
\begin{split}
& \sum_{t=-\infty} ^{\infty} (V^* U_0^{t+1} E_{U_0} (I) f, U^{-\tau} V^* U_0^{t+\tau +1} E_{U_0} (I)g ) _{\mathcal{H}} \\
&= \int_I \sum_{t=-\infty}^{\infty } ( \mathcal{F}_0 (\theta ) VU^{\tau} V^* U_0 ^{t+1} E_{U_0} (I) f, e^{i(t+\tau +1 )} \mathcal{F}_0 (\theta )g) _{{\bf h} (\theta )} d\theta .
\end{split}
\end{gather*}
Here we use the formula (\ref{S5_eq_WO2}) again, inserting $e^{-\epsilon |t+1|} $ and letting $\epsilon \downarrow 0$.
Thus the above sum is equal to 
$$
2\pi \int_I e^{-i\tau \theta} ( \mathcal{F}_0 (\theta ) VU^{\tau} V^* \mathcal{F}_0 (\theta )^* \mathcal{F}_0 (\theta ) f, \mathcal{F}_0 ( \theta )g) _{{\bf h} (\theta )} d\theta ,
$$
and then
$$
{\bf S}_2 = 2\pi \int_I \sum _{\tau =0}^{\infty} ( \mathcal{F}_0 (\theta ) V(e^{-i\tau \theta } U^{\tau}) V^* \mathcal{F}_0 (\theta )^* \mathcal{F}_0 (\theta ) f, \mathcal{F}_0 ( \theta )g) _{{\bf h} (\theta )} d\theta .
$$
The formula (\ref{S5_eq_WO1}) implies $ \sum_{\tau =0}^{\infty} e^{-i\tau \theta} U^{\tau} = -e^{i\theta} R(\theta -i0 )$.
Then
\begin{equation}
{\bf S}_2 = - 2\pi \int_I e^{i\theta} ( \mathcal{F}_0 (\theta ) V R(\theta -i0 ) V^* \mathcal{F}_0 (\theta ) ^* \mathcal{F}_0 (\theta )f, \mathcal{F}_0 (\theta )g ) _{{\bf h} (\theta )} d\theta .
\label{S5_eq_Smatrix12}
\end{equation}
The equalities (\ref{S5_eq_Smatrix11}) and (\ref{S5_eq_Smatrix12}) implies (\ref{S5_eq_Smatrix_ampli}).
\qed

\subsection{Singularity expansion of generalized eigenfunction}

We have arrived at the last topic of the present paper.
We consider a singularity expansion of the generalized eigenfunction $ \mathcal{F}_+ (\theta )^* \phi $ for $ \phi \in {\bf h} (\theta )$.
The discussion in this subsection will imply that the S-matrix $ \widehat{{\bf S} } (\theta )$ is an analogue of the physical S-matrix as a mapping from the incident wave to the scattered wave in view of the scattering theory of Schr\"{o}dinger equations.

Note that 
$$
R(\theta +i0)^* = -e^{i\theta} U R( \theta -i0) .
$$
For $ \phi \in {\bf h} (\theta )$ with $ \theta \in J_{\gamma} \setminus J_{\gamma , \mathcal{T}} $, we have 
\begin{gather*}
\begin{split}
\mathcal{F}_+ (\theta )^* \phi &= \mathcal{F}_0 (\theta )^* \phi - R( \theta +i0 )^* V^* \mathcal{F}_0 (\theta )^* \phi \\
&= \mathcal{F}_0 (\theta )^* \phi + e^{i\theta} U R( \theta -i0 ) V^* \mathcal{F}_0 (\theta )^* \phi \\
&=\mathcal{F}_0 (\theta )^* \phi + e^{i\theta} \left( e^{i\theta} R(\theta -i0) V^* \mathcal{F}_0 (\theta )^* \phi + V^* \mathcal{F}_0 (\theta )^* \phi \right) .
\end{split}
\end{gather*}
Then the generalized eigenfunction $ \mathcal{F}_+ (\theta )^* \phi $ satisfies
$$
\mathcal{F}_+ (\theta )^* \phi = \mathcal{F}_0 (\theta )^* \phi + e^{2i\theta} R(\theta -i0) V^* \mathcal{F}_0 (\theta )^* \phi ,
$$
up to rapidly decreasing terms.
It follows from Lemma \ref{S3_lem_parseval} and Theorem \ref{S3_thm_rangeF0*} that there exists $f\in \mathcal{B} ({\bf Z})$ such that $ \phi = \mathcal{F}_0 (\theta )f$ and
$$
R_0 (\theta +i0 )f -R_0 (\theta -i0)f = 2\pi e^{-i\theta} \mathcal{F}_0 (\theta )^* \mathcal{F}_0 (\theta )f= 2\pi e^{-i\theta} \mathcal{F}_0 (\theta )^* \phi .
$$
Then we have 
\begin{gather}
\begin{split}
&\mathcal{F}_+ (\theta )^* \phi \\
&= \frac{e^{i\theta}}{2\pi} \left( R_0 (\theta +i0 ) f - R_0 (\theta -i0 )f +2\pi e^{i\theta} R (\theta -i0) V^* \mathcal{F}_0 (\theta )^* \phi \right) \\
&=\frac{e^{i\theta}}{2\pi} \left( R_0 (\theta +i0 )f - R_0 (\theta -i0) \left( f - 2\pi e^{i\theta} (V^* - VR(\theta -i0) V^*) \mathcal{F}_0 (\theta )^* \phi    \right)   \right) ,
\end{split}
\label{S5_eq_F+*formula}
\end{gather} 
up to rapidly decreasing terms.
Now let us define
\begin{gather}
\phi^{in}: = \mathcal{F}_0 (\theta )f = \phi , \label{S5_def_incoming} \\
\phi ^{out} := \mathcal{F}_0 (\theta )\left( f - 2\pi e^{i\theta} (V^* - VR(\theta -i0) V^*) \mathcal{F}_0 (\theta )^* \phi    \right) =\widehat{{\bf S}} ( \theta ) \phi ^{in} . \label{S5_def_outgoing}
\end{gather}
Summarizing the above computation, we obtain the following fact.

\begin{lemma}
Let $ \phi ^{in} , \phi ^{out} \in {\bf h} (\theta )$ are defined by (\ref{S5_def_incoming}) and (\ref{S5_def_outgoing}).
Then we have $ \phi^{out} = \widehat{{\bf S}} (\theta )\phi^{in} $.

\label{S5_lem_Smatrixpattern}
\end{lemma}

At the end of this section, we derive a characterization of the generalized eigenfunction of the equation $ (U-e^{i\theta} )u=0$ for $ \theta \in J_{\gamma} \setminus J_{\gamma , \mathcal{T}} $.

\begin{lemma}
For $ \theta \in J_{\gamma } \setminus J_{\gamma , \mathcal{T} }$ and $ \phi \in {\bf h} (\theta )$, there exist some constants $\kappa_2 (\theta ) > \kappa_1 (\theta )>0$ such that 
$$
\kappa_1 (\theta ) \| \phi \|_{{\bf h} (\theta )} \leq \| \mathcal{F}_+ (\theta )^* \phi \| _{\mathcal{B}^* ({\bf Z})} \leq \kappa_2 (\theta ) \| \phi \| _{{\bf h} (\theta )} .
$$

\label{S5_lem_F+closedrange}
\end{lemma}

Proof.
In view of Lemma \ref{S4_lem_stone}, we have only to prove the estimate from below.
Recalling $ \mathcal{F}_0 (\theta ) \mathcal{B} ({\bf Z}) = {\bf h} (\theta )$, we take $ f^{in} , f^{out} \in \mathcal{B} ({\bf Z})$ such that $ \mathcal{F}_0 (\theta ) f^{in} = \phi ^{in} := \phi $ and $ \mathcal{F}_0 (\theta ) f^{out} = \phi^{out} := \widehat{{\bf S} } (\theta ) \phi ^{in} $.
In view of (\ref{S5_eq_F+*formula}), we have 
$$
\mathcal{F}_+ (\theta )^* \phi  = \frac{e^{i\theta}}{2\pi} \left( R_0 (\theta +i0 ) f^{in} - R_0 (\theta -i0 ) f^{out}  \right) .
$$
Let $\widehat{u} = \mathcal{U} \mathcal{F}_+ (\theta )^* \phi $ and $ \xi (\theta )\in M(\theta )$.
We take $ \chi \in C^{\infty} ({\bf T})$ such that $\chi (\xi (\theta )) =1$ with small support.
We put $\widehat{u} _{j} = \chi P_j \widehat{u} $. 
Then we introduce the change of variable $ \xi \mapsto \eta $ as in \S 3 in a neighborhood of $ \xi (\theta ) \in M(\theta )$.

In the following, we shall prove for $ \theta \in J_{\gamma,1} \setminus J_{\gamma , \mathcal{T} }$.
For $ \theta \in J_{\gamma ,2} \setminus J_{\gamma , \mathcal{T}} $, the proof is similar.
It follows from (\ref{S3_eq_LAP12}) and (\ref{S3_eq_LAP13}) that 
$$
\widetilde{u}_1 (s)= \frac{1}{2\pi} \left( \int_{s+1}^{\infty} \widetilde{f}_1^{in} (t) dt + \int_{-\infty}^{s+1} \widetilde{f}_1^{out} (t) dt  \right) ,
$$
up to terms in $ \mathcal{B}_0^* ({\bf R}) $ where $\widehat{f}_1^{in} = \chi P_1 \widehat{f}^{in}$ and $ \widehat{f}_1^{out} = \chi P_1 \widehat{f}^{out} $.
Moreover, Lemma \ref{S3_lem_LAP1} implies that 
$$
\widetilde{u}_1 (s)= \frac{1}{2\pi} \left( H(-s-1) \int_{-\infty}^{\infty} \widetilde{f}_1^{in} (t) dt + H(s+1) \int _{-\infty}^{\infty} \widetilde{f}_1^{out} (t) dt \right) ,
$$
up to terms in $\mathcal{B}_0^* ({\bf R})$.
Thus we have
$$
\limsup _{R\to \infty} \frac{1}{R} \int_{|s|<R} |\widetilde{u}_1 (s)|^2 _{{\bf C}^2} ds \geq \kappa \left( \left| \int_{-\infty}^{\infty} \widetilde{f}_1^{in} (t)dt \right|^2 _{{\bf C}^2} +\left| \int_{-\infty}^{\infty} \widetilde{f}_1^{out} (t)dt \right|^2 _{{\bf C}^2}   \right) ,
$$
for a constant $\kappa >0$, noting 
$$
\mathrm{Re} \left( \int_{s+1}^{\infty} \widetilde{f}_1^{in} (t) dt \int_{-\infty}^{s+1} \overline{\widetilde{f}_1^{out} (t)} dt \right) \in \mathcal{B}_0 ^* ({\bf R}). 
$$
From the equality 
$$
\int_{-\infty}^{\infty} \widetilde{f}_1^{in} (t)dt = \sqrt{2\pi} \widehat{f}_1 ^{in} (\eta ) \big| _{\eta =0} , \quad \int_{-\infty}^{\infty} \widetilde{f}_1^{out} (t)dt = \sqrt{2\pi} \widehat{f}_1 ^{out} (\eta ) \big|_{\eta =0} ,
$$
and the inequality (\ref{S3_eq_ineqB*limsup}), we obtain
\begin{equation}
\sup _{R>1} \frac{1}{R} \int _{|s|<R} |\widetilde{u}_1 (s)|^2 _{{\bf C}^2} ds \geq \kappa \left( |\widehat{f} _1 ^{in} (0)|_{{\bf C}^2}^2 + |\widehat{f} _1 ^{out} (0)|_{{\bf C}^2}^2 \right) ,
\label{S5_eq_ineqB*}
\end{equation}
for a constant $\kappa >0$.
The inequality (\ref{S5_eq_ineqB*}) means in the variable $\xi$ that
\begin{gather}
\begin{split}
\| \widehat{u} \|^2_{\mathcal{B}^* ({\bf T})} & \geq \kappa \left( \| \mathcal{F}_0 (\theta )f^{in} \| _{{\bf h}(\theta )} ^2 + \| \mathcal{F}_0 (\theta )f^{out} \| _{{\bf h}(\theta )} ^2 \right) \\
&= \kappa \left( \| \phi^{in} \| _{{\bf h}(\theta )} ^2 + \| \phi^{out} \| _{{\bf h}(\theta )} ^2 \right)  .
\end{split}
\label{S5_eq_ineqB*2}
\end{gather}
Since $\widehat{{\bf S}} (\theta )$ is unitary on ${\bf h} (\theta )$ and $\phi^{out} = \widehat{{\bf S}} (\theta )\phi ^{in}$, the above inequality proves the lemma.
\qed

\begin{theorem}
Let $ \theta \in J_{\gamma} \setminus J_{\gamma , \mathcal{T}} $.
Then we have $\mathcal{F}_+ (\theta ) \mathcal{B} ({\bf Z} )={\bf h } (\theta ) $ and $ \{ u\in \mathcal{B}^* ({\bf Z}) \ ; \ (U-e^{i\theta} )u=0 \} = \mathcal{F}_+ (\theta )^* {\bf h} (\theta )$.

\label{S5_thm_rangeF+}
\end{theorem}

Proof.
In view of Lemma \ref{S5_lem_F+closedrange}, the range of $\mathcal{F}_+ (\theta )^* $ is closed. 
Thus we use Theorem \ref{S3_thm_Banach} again, letting $X_1 = \mathcal{B} ({\bf Z})$, $X_2 = {\bf h} (\theta )$ and $T= \mathcal{F}_+ (\theta )$.
We have only to prove $(u,f)=0$ when $u\in \mathcal{B} ^* ({\bf Z})$, $f\in \mathcal{B} ({\bf Z})$, $(U-e^{i\theta} )u=0$ and $\mathcal{F}_{+} (\theta )f=0$. 

We put $v= R(\theta - i0)^* f$.
Note that $v\in \mathcal{B}^* ({\bf Z})$ satisfies $ (U^* -e^{-i\theta} )v=f$ i.e. $ (U-e^{i\theta} )v= -e^{i\theta} Uf$. 
Thus $v= -R(\theta + i0)(Uf) $.
By the assumption $ \mathcal{F}_{+} (\theta ) f=0$, we have $ \mathcal{F} _{+} (\theta ) (Uf)= e^{i\theta} \mathcal{F}_{+} (\theta ) f=0$.
Thus we have $v\in \mathcal{B}_0^* ({\bf Z})$ from the resolvent equation and the assertion (4) of Theorem \ref{S3_thm_LAPlattice}. 

Let us compute $ (u,f)=(u, (U^* -e^{-i\theta} )v )$. 
Take $ \psi \in C_0^{\infty} ({\bf R})$ such that $\psi (t)=1 $ for $t<1$, $\psi (t)=0$ for $ t>2$.
Then we have for large $R>1$
\begin{equation}
(u, \psi (|X|/R) (U^* -e^{-i\theta} )v) = ([U, \psi (|X|/R) ] u,v) ,
\label{S4_eq_characterize11}
\end{equation}
by using $ (U-e^{i\theta} )u=0$.
We note that 
\begin{gather*}
\begin{split}
[U, \psi (|X|/R)] &= [ U_0 , \psi (|X|/R )] + [V, \psi (|X|/R)] \\
&= \Psi (x) (SC_0 + S(C(x)-C_0 )) ,
\end{split}
\end{gather*}
where
\begin{gather*}
\Psi (x)= \left[ \begin{array}{cc}
\psi (|X+1|/R) -\psi (|X|/R) & 0 \\ 0 & \psi (|X-1|/R) - \psi (|X|/R) \end{array} \right]  .
\end{gather*}
Let $ F_{x,\pm} (t)= \psi (|x \pm t|/R) - \psi (|x|/R)$ for $x \in {\bf Z} $ and $ t\in [0,1]$.   
In view of Taylor's theorem, we have 
\begin{gather*}
\begin{split}
F_{x,\pm} (1) &= \pm \frac{\mathrm{sgn}  (x) }{R} \psi' (|x|/R) + \frac{1}{R^2} \int_0^1 (1-s) \psi '' (|x|/R) ds  \\
&=\pm \frac{\mathrm{sgn}  (x) }{R} \psi' (|x|/R) +o(R^{-1}) ,
\end{split}
\end{gather*}
as $R\to \infty $ for $x \not= 0$.
This implies 
$$
\Psi (x)= \frac{\mathrm{sgn} (x)}{R} \psi ' (|X|/R) \left[\begin{array}{cc} 1 & 0 \\ 0 & -1 \end{array} \right] +o(R^{-1} ) ,
$$
for large $|x| $ and $R$. 
The equality (\ref{S4_eq_characterize11}) can be estimated by
$$
\kappa \left( \frac{1}{R} \sum _{|x|<\kappa R} |u(x)|^2 \right)^{1/2} \left( \frac{1}{ R} \sum _{|x|<\kappa R} |v(x)|^2  \right) ^{1/2} ,
$$
for some constants $ \kappa >0$.
Tending $R\to \infty$ on this equality, we obtain $ (u,f)=0$.
\qed

\appendix
\section{Some complex contour integrations}

Let $\epsilon >0$ be small.
We derive formulas for the functions
\begin{equation}
\widetilde{F} _{  \pm, \epsilon } (s)= \int _{-\pi}^{\pi} \frac{ e^{-is\eta}}{e^{i\eta} -1\pm \epsilon } d\eta , \quad  s\in {\bf R} .
\label{A_def_Fpme}
\end{equation}
We put
\begin{gather}
I_+ (s)= H(-s-1) \frac{\sin (\pi s)}{\pi} \int ^{\infty}_0 \frac{e^{s\tau}}{e^{-\tau} +1} d\tau , \label{A_def_I+} \\
I_- (s)= H(s+1) \frac{\sin (\pi s)}{\pi} \int _{-\infty}^0 \frac{e^{s\tau}}{e^{-\tau} +1} d\tau ,\label{A_def_I-}
\end{gather}
where $H(s)$ is the Heaviside function.

\begin{lemma}
We have $ I_{\pm}  \in L^{\infty} ({\bf R}) \cap L^2 ({\bf R})$.
\label{A_lem_L2}
\end{lemma}

Proof.
For $s+1<0$, we have 
$$
| I_+ (s) |  \leq \frac{| \sin (\pi s) |}{\pi} \int ^{\infty}_0 \frac{e^{(s+1)\tau}}{1+e^{\tau} } d\tau 
\leq \frac{| \sin (\pi ( s+1)) |}{\pi |s+1|}   .
$$
Thus $|I_+ (s)| = O(|s|^{-1} ) $ as $|s| \to \infty $, and $|I_+ (s)|$ is bounded in a neighborhood of $s+1 =0$. 
This implies $I_+ \in L^2 ({\bf R})$. 
For $ I_-$, the proof is similar.
\qed

\begin{lemma}
We have
$$
\lim_{\epsilon \downarrow 0} \widetilde{F}_{\pm ,\epsilon} (s)=  2\pi   \left( \pm  H(\mp (s+1)) + I_+ (s)   - I_- (s ) \right) .
$$
\label{A_lem_fourier}
\end{lemma}

Proof.
For the proof, we use the contour integration on the complex plane.
Let us prove the formula for $\widetilde{F}_{+,\epsilon} $.
If $s+1<0$, we take a path $C= C_1 + C_2 + C_3 + C_4$ such that 
$$
C_1 =\{ z=\tau \ ; \ \tau : -\pi \to \pi \} , \quad C_2 =\{ z= \pi +i\tau \ ; \ \tau :0\to \rho \} , 
$$
$$
C_3 =\{ z=\tau +i\rho \ ; \ \tau : -\pi \to \pi \} ,\quad C_4 =\{ z= -\pi + i\tau \ ; \ \tau : \rho \to 0 \} ,
$$
with large $\rho >0 $.
Letting 
$$
f_{\epsilon} (z)= \frac{e^{-isz}}{e^{iz} -1+\epsilon }  , \quad z\in {\bf C} ,
$$
the function $f_{\epsilon} $ has a simple pole at $z=-i \log (1-\epsilon )$ and its residue is $-i (1-\epsilon )^{-s-1} $.
The residue theorem implies
\begin{equation}
\int_C \frac{e^{-isz}}{e^{iz} -1+\epsilon} dz = 2\pi (1-\epsilon )^{-s-1} .
\label{A_eq_residue}
\end{equation}
On $C_2 $ and $C_4 $, we have
$$
\int _{C_2} \frac{e^{-isz}}{e^{iz} -1+\epsilon } dz = -i e^{-is\pi}  \int _0^{\rho} \frac{ e^{s\tau} }{e^{-\tau} +1-\epsilon } d\tau ,
$$
and
$$
\int _{C_4} \frac{e^{-isz}}{e^{iz} -1+\epsilon } dz = i e^{is\pi}  \int _0^{\rho} \frac{ e^{s\tau} }{e^{-\tau} +1-\epsilon } d\tau .
$$
Then we have
\begin{equation}
\int _{C_2} + \int _{C_4} \frac{e^{-isz}}{e^{iz} -1+\epsilon } dz = -2 \sin (\pi s ) \int _0^{\rho} \frac{e^{s\tau}}{e^{-\tau} +1-\epsilon } d\tau .
\label{A_eq_C2C4}
\end{equation}
On $C_3$, we obtain
\begin{equation}
\left| \int _{C_3}  \frac{e^{-isz}}{e^{iz} -1+\epsilon } dz \right| \leq \frac{2\pi e^{(s+1) \rho }}{e^{\rho} (1-\epsilon ) -1 } .
\label{A_eq_C3}
\end{equation}
Plugging (\ref{A_eq_residue})-(\ref{A_eq_C3}) and tending $\rho \to \infty $, we have
\begin{equation}
\int_{-\pi}^{\pi} \frac{e^{-is\eta} }{e^{i\eta} -1+\epsilon } d\eta = 2\pi (1-\epsilon )^{-s-1} +2\sin (\pi s) \int _0^{\infty} \frac{e^{s\tau}}{e^{-\tau} +1-\epsilon } d\tau , 
\label{A_eq_F+e}
\end{equation}
for $s+1<0 $.

If $s+1>0$, we take a path $C' = C' _1 + C' _2 + C' _3 + C' _4 $ such that 
$$
C' _1 = \{ z=\tau \ ; \ \tau : \pi \to -\pi \} , \quad C'_2 =\{ z= -\pi + i\tau \ ; \ \tau : 0\to -\rho \} ,
$$
$$
C' _3 =\{ z= \tau -i\rho \ ; \ \tau : -\pi \to \pi \} , \quad C' _4 = \{ z=\pi +i\tau \ ; \ \tau : -\rho \to 0 \} .
$$
Cauchy's integration theorem implies 
\begin{equation}
\int _{C'} \frac{e^{-isz}}{e^{iz} -1+\epsilon } dz =0.
\label{A_eq_cauchy}
\end{equation}
As has been in the above argument, we also have
\begin{equation}
\int _{C'_2} + \int _{C' _4} \frac{ e^{isz} }{e^{iz} -1+\epsilon } dz = -2 \sin (\pi s ) \int _{-\rho}^0 \frac{e^{s\tau} }{e^{-\tau} +1 -\epsilon } d\tau ,
\label{A_eq_C2C4prime}
\end{equation}
and
\begin{equation}
\left| \int _{C'_3} \frac{e^{-isz}}{e^{iz} -1+\epsilon } dz \right| \leq \frac{2\pi e^{-(s+1)\rho } }{1-e^{-\rho} (1-\epsilon )} .
\label{A_eq_C3prime}
\end{equation}
Tending $\rho \to \infty $, the equations (\ref{A_eq_cauchy})-(\ref{A_eq_C3prime}) imply 
\begin{equation}
\int _{-\pi}^{\pi} \frac{ e^{-is\eta} }{e^{i\eta} -1+\epsilon } d\eta = -2 \sin (\pi s ) \int _{-\infty}^0 \frac{e^{s\tau}}{e^{-\tau} +1 -\epsilon } d\tau .
\label{A_eq_F-e}
\end{equation}
The formulas (\ref{A_eq_F+e}) and (\ref{A_eq_F-e}) prove the lemma for $\widetilde{F} _{+,\epsilon }$.
For $\widetilde{F}_{-,\epsilon} $, the proof is similar.
\qed

\medskip

Let us turn to the integration
\begin{equation}
S_{\pm , \epsilon} (\omega )= \int _{\theta_1}^{\theta_2}  \frac{d \theta }{e^{i(  \omega -\theta )} -1 \pm \epsilon}  , \quad \epsilon >0 ,\quad 0< \theta_1 < \theta_2 < 2\pi .
\label{A_eq_stone}
\end{equation}

\begin{lemma}
We have 
\begin{gather*}
\lim _{\epsilon \downarrow 0} ( S_{+,\epsilon} (\omega ) - S_{-,\epsilon } (\omega )) = \left\{ 
\begin{split}
2\pi &, \quad \omega \in (\theta _1 , \theta _2 ) , \\
0 & , \quad \omega \in (-\infty , \theta_1 ) \cup ( \theta_2 , \infty ) .
\end{split} \right.
\end{gather*}
\label{A_lem_stoneformula}
\end{lemma}

Proof.
We rewrite $S_{\pm , \epsilon} ( \omega )$ as 
$$
S_{\pm , \epsilon} ( \omega )=  \int _{\omega -\theta_2}^{\omega -\theta_1} \frac{ d\eta }{ e^{i\eta} -1 \pm \epsilon }  .
$$
We use the contour integration on the complex plane again.
We take a path $C=C_1 + C_2 + C_3 + C_4$ such that 
$$
C_1 =\{ z=\tau \ ; \ \tau : \theta -\theta_1 \to \theta -\theta_2 \} ,\quad C_2 =\{ z= \theta -\theta_2 +i\tau \ ; \ \tau : 0\to -\rho \} ,
$$
$$
C_3 =\{ z=\tau -i\rho \ ; \ \tau : \theta -\theta_2 \to \theta -\theta_1 \} , \quad C_4 =\{ z= \theta -\theta_1 +i\tau \ ; \ \tau : -\rho \to 0\} ,
$$
for sufficiently large $\rho >0$. 
Thus we compute as in the proof of Lemma \ref{A_lem_fourier}.
We have
$$
S_{+,\epsilon} (\omega ) = I_{+,1, \epsilon} (\omega ) - I_{+,2 ,\epsilon } (\omega )  ,
$$
and
\begin{gather*}
S_{-,\epsilon} (\omega )  = \left\{
\begin{split}
 -\frac{2\pi }{1+\epsilon}  + I_{+,2,\epsilon} (\omega ) - I_{-,2,\epsilon} (\omega ) & , \quad \omega \in (\theta_1 , \theta_2 ) ,\\
 I_{+,2,\epsilon} (\omega ) - I_{-,2,\epsilon} (\omega ) & , \quad \omega \in (-\infty , \theta_1 )\cup  ( \theta_2 , \infty ) .
\end{split}
\right.
\end{gather*}
where 
$$
 I_{j,\pm ,\epsilon } (\omega )= i\int _{-\infty}^0 \frac{d\tau}{e^{i(\omega - \theta_j )} e^{-\tau} -1 \pm \epsilon } . 
$$
Then we have proven the lemma.
\qed

\end{document}